\documentclass[11pt]{article}

\usepackage[all]{xy}

\usepackage{mathrsfs} %
\usepackage{bbm}
\usepackage{enumerate} %

\usepackage{amsfonts}
\usepackage{amssymb}

\renewcommand{\caption}[1]{\centerline{{\bf Table \thetable}. {#1}}}

\newtheorem{theorem}{Theorem}%
\newtheorem{prop} [theorem]{Proposition}
\newtheorem{cor}  [theorem]{Corollary}
\newtheorem{example}[theorem]{Example}
\newtheorem{lemma}[theorem]{Lemma}
\newtheorem{defn} [theorem]{Definition}
\newtheorem{rem}  [theorem]{Remark}

\newcommand{\com}{{\mathbb C}}

\newcommand{\zet}{{\mathbb Z}}

\newcommand{\bin}{{\mathbb B}}
\newcommand{\quat}{{\mathbb H}}
\newcommand{\oct}{{\mathbb O}}
\newcommand{\binz}{{\mathbb \bin_\zet}}
\newcommand{\quatz}{{\mathbb \quat_\zet}}
\newcommand{\octz}{{\mathbb \oct_\zet}}

\newcommand{\one}{{\mathbbm 1}}

\newcommand{\calg}{{C}}

\newcommand{\J}{{\mathfrak J}}
\newcommand{\aJ}{\underline{\mathfrak J}}
\newcommand{\m}{{\mathfrak M}}          %
\newcommand{\jz}{{\mathfrak \J_\zet}}    %
\newcommand{\mz}{{\mathfrak \m_\zet}}    %
\newcommand{\mj}{\ensuremath {\m(\J)}}   %

\newcommand{\cart}{{\mathfrak h}}
\newcommand{\lie}{{\mathfrak g}}
\newcommand{\ml}{{\mathfrak m}}

\newcommand{\adlie}{{\rm ad}}

\newcommand{\dvap}{{d\,'_2}}

\newcommand{\alg}[1]{\underline{#1}}

\newcommand{\npg}[1]{{{\rm NP}(#1)}}     %
\newcommand{\str}[1]{{{\rm Str}(#1)}}
\newcommand{\astr}[1]{{\alg{\rm Str}(#1)}}
\newcommand{\anpg}[1]{{\alg{\rm NP}(#1)}}

\newcommand{\g}[1]{{\rm Inv}({#1})}     %
\newcommand{\gr}{{\g{\m}}}     %
\newcommand{\ag}[1]{{\alg{\rm Inv}(#1)}}     %

\newcommand{\grz}{{\g{\mz}}}   %

\newcommand{\GL}{{\rm GL}}
\newcommand{\SL}{{\rm SL}}
\newcommand{\SO}{{\rm SO}}

\newcommand{\cross}{{\times}} %

\newcommand{\h}[1]{{\cal H}_3(#1)}

\newcommand{\mto}{\mapsto}
\newcommand{\ep}{$\blacksquare$}
\newcommand{\ee}{$\square$}
\newcommand{\tensor}{\otimes}

\newcommand{\bfa}{{\bf a}}
\newcommand{\bfb}{{\bf b}}
\newcommand{\bfmu}{\mbox{\boldmath $\mu$}}

\newcommand{\mnorm}{N}  %
\newcommand{\ncomp}{{\bf n}}  %
\newcommand{\tcomp}{{\bf t}}  %
\newcommand{\tr}{{\rm Tr}} %
\newcommand{\btr}{} %

\newcommand{\np}{norm-preserving}

\newcommand{\rank}{{\rm rank\,}}

\newcommand{\jprod}{\bullet}

\newcommand{\cc}{{}^\circ}

\newcommand{\conj}[1]
{ \overline{#1} }
\newcommand{\defeq}{ \stackrel{\scriptscriptstyle \rm def}{=} }

\newcommand{\mel}[4]   %
{ \maketwo   {#1}{#3}{#4}{#2}   }

\newcommand{\maketwo}[4]
{
\left(
\begin{array}{cc}
   #1 & {#2} \\
   #3 & #4
\end{array}
\right)
}

\newcommand{\maketwoblock}[4]
{
\left(
\begin{array}{c|c}

   #1 & #2 \\
\hline
   #3 & #4
\end{array}
\right)
}

\newcommand{\makethree}[9]
{ \left(
\begin{array}{ccc}

   #1        & #2        & #3 \\
   #4        & #5        & #6 \\
   #7        & #8        & #9
\end{array}
\right)
}

\newcommand{\makeherm}[6]
{ \left(
\begin{array}{ccc}

   #1        & #6        & \conj{#5} \\
   \conj{#6} & #2        & #4        \\
   #5        & \conj{#4} & #3
\end{array}
\right)
}
\newcommand{\bx}{{x}}
\newcommand{\id}{{\rm \,Id}}

\newcommand{\F}{F}                      %

\newcommand{\tf}[1]{T({#1})}             %

\newcommand{\mod}{{\rm mod}\,}

\newcommand{\mult}{\times} %

\renewcommand{\gcd}{\ensuremath{g.c.d.\,}}

\newcommand{\diag}[3]{{{\rm diag}\{#1, #2, #3\}}}

\newcommand{\todoz}[1]{}
\newcommand{\myskip}[1]{}

\newcommand{\End}{{\rm End}}
\newcommand{\mat}[4]{
    \left[
    \begin{array}{cc}
       #1  &  #2 \\
       #3  &  #4
    \end{array}
    \right]
}

\newcommand{\matthree}[9]{
    \left[
    \begin{array}{ccc}
       #1  &  #2  & #3 \\
       #4  &  #5  & #6 \\
       #7  &  #8  & #9
    \end{array}
    \right]
}

\begin{document}
\setcounter{table}{1}

\begin{center}
{
\Large\bf

Jordan algebras, exceptional groups,\\[2mm]
and higher composition laws
}

\vspace{\baselineskip}

{\large Sergei Krutelevich}

\vspace{\baselineskip}

\end{center}

\begin{abstract}

We consider an integral version of the Freudenthal construction relating
Jordan algebras and exceptional algebraic groups. We show how this
construction is related to higher composition laws of M.~Bhargava in
number theory \cite{bh1}.

We propose an algorithmic approach to studying orbit spaces of
groups underlying higher composition laws.
Using this method we discover two new
examples of spaces sharing similar properties, and indicate
several more examples of spaces where such composition laws may be
introduced.
\end{abstract}

\tableofcontents

\section{Introduction}
\subsection{Integral representations of exceptional groups}
\label{i:fc}

It is a well-known fact in number theory that there is a one-to-one
correspondence between
the set of $\SL_2(\zet)$-equivalence classes of integral binary quadratic
forms and the set of (narrow) ideal classes in quadratic rings. About two
hundred years ago Gauss discovered the law of composition of binary
quadratic forms, which turns
the set of equivalence classes of primitive forms of a given discriminant
$D$ into a group isomorphic to the ideal class group of the quadratic ring
of discriminant $D$. This correspondence is a very important tool
for doing computations in the ideal class group.

A few years ago M.~Bhargava discovered several more examples of the same
kind, which he referred to as higher composition laws. More precisely, he
showed that there are other examples of linear groups $G_\zet$
and their integral representations $V_\zet$ such that $G_\zet$-orbits in $V_\zet$
are in one-to-one correspondence with ideal classes in the rings of
integers in number fields, see \cite{bh-phd, bhau, bh1}.
M.~Bhargava also noted a surprising connection between spaces underlying higher
composition and exceptional Lie groups.

In the present paper we investigate this connection with exceptional
groups. We show that there is a natural construction which assigns
a cubic Jordan algebra to every space underlying higher composition laws
associated with quadratic rings. We study the appropriate orbit spaces using an
algorithmic procedure reminiscent of the Gaussian elimination algorithm
for the usual integral matrices. This approach allows us to provide two
new examples of spaces sharing similar properties, and indicate
several more examples of spaces where such composition laws may be
introduced.

\medskip
Our interest in the study of integral representations of exceptional groups was
originally motivated by a question on a standard form of a charge vector in a certain
physical model, which appeared in a paper on BPS black holes in string theory by
H.~Maldacena, G.~Moore, A.~Strominger \cite{mms}.
This question is equivalent to the question on a normal form of a vector
in the $27$-dimensional integral representation of the split group of type
$E_6$.
Such a representation may be constructed via the exceptional Jordan
algebra of $3\times 3$ Hermitian matrices over the algebra of octonions.
The answer, generalizing the classical results on the Smith normal
form for regular integral matrices, was given in our earlier paper
\cite{k}, see also Theorem~$\ref{np-orbits}$ below.

It was suggested by G.~Moore that another interesting space to look
at in this context is the $56$-dimensional integral representation of the
split group of type $E_7$. This representation may be constructed using the
$27$-dimensional exceptional Jordan algebra. This procedure is known
as the Freudenthal construction \cite{b, fr14},
and it can be applied to other cubic Jordan algebras.
Given a cubic Jordan algebra~$\J$, the Freudenthal  construction
produces a semisimple algebraic group and an irreducible module over it.
We will denote the group obtained by $\gr$ and the module by $\mj$
(see Subsection~\ref{ss:fc} for details).
Moreover, this module is equipped with a quartic form and a symplectic
form, which are invariant under the action of the group.
When $\J$ is the $27$-dimensional exceptional Jordan
algebra, the group $\gr$ is $E_7$ and $\mj$ is its $56$-dimensional
representation.

A natural
way to construct a Jordan algebra with a cubic form is to consider the
space of $3\times 3$-Hermitian matrices over a composition algebra
over a field $F$, the cubic form being the determinant of such matrices
(see Example~$\ref{ex:jcubic}$ for details). When $F$ is algebraically
closed, there are four composition algebras, and they produce four cubic
Jordan algebras of dimension $6, 9, 15, 27$.

The application of the Freudenthal construction yields a module of
dimension $14, 20, 32, 56$ for a certain simple algebraic group
of type $C_3, A_5, D_6, E_7$.
This module is in fact a prehomogeneous vector space for this
algebraic group, and it is natural to study orbits under the action of
the group.
In the case $F=\com$, the classification of orbits of the one-dimensional
subspaces\footnote{We avoid
using the term ``projective" here, since it is used in completely
different sense in the rest of the paper.}
arising this way was obtained by J.-L.~Clerc~\cite{clerc}.

Lie algebras of the algebraic groups of type $C_3, A_5, D_6, E_7$ appear
in the third row of the Freudenthal-Tits magic square. The groups,
produced by the Freudenthal construction, also
appear in the last row of the
``magic triangle" of subgroups associated to the exceptional
series of P.~Deligne, see \cite{dg}. In addition, the representation $\m$
of the group $\gr$ is the {\em preferred} representation of this group
in the sense of \cite{dg}.

In the present paper we consider an integral version of the Freudenthal
construction.
In this case the quartic form (the {\em norm}) takes on integer values.
We study integral forms of the split groups of type $A_5, D_6, E_7$
determined by the Freudenthal construction, and
we obtain the following result
on the structure of integral orbits in the $\zet$-modules of
dimension $20, 32, 56$.
\medskip

\hspace{-\parindent}{\bf Theorem~$\ref{thm:main}$}
{\em
Let $(G_\zet, \mz)$ be one the following pairs
$$
\Bigl(\SL_6(\zet), \wedge^3(\zet^6)\Bigr),\quad
\Bigl(D_6(\zet), \mbox{\rm half-spin}_\zet\Bigr),\quad
\Bigl(E_7(\zet), V(\omega_7)_\zet\Bigr),
$$
Then
\begin{itemize}

\item
The $G_\zet$-invariant quartic form (the norm) on the module $\mz$
has values congruent to $0$ or $1\, (\!\mod 4)$.

\item
Let $n$ be an integer $\equiv 0$ or $1\, (\!\mod 4)$.
The group $G_\zet$ acts
transitively on the set of {\em projective} elements of norm $n$.

\item
If $n$ is a fundamental discriminant\footnote{
An integer $n$ is called a fundamental discriminant if $n$ is squarefree
and $\equiv 1 (\mod 4)$
or $n=4k$, where $k$ is a squarefree integer that is $\equiv 2$ or $3 (\mod 4)$.
The result stated in the theorem applies in the case $n=1$.
},
then every element of norm $n$ is projective, and hence in this case
$G_\zet$ acts transitively on the set of elements of norm $n$.
\end{itemize}
}

The assertion of this theorem for $\SL_6(\zet)$-orbits
on $\wedge^3(\zet^6)$ (when $n\ne 0$) was proved by M.~Bhargava in
\cite[Theorem 7]{bh1}, using the correspondence with (narrow) ideal
classes in quadratic orders.

Our approach, based on the Freudenthal construction, allows us to treat
all values of the norm (including $n=0$) uniformly.
Our statement of Theorem~$\ref{thm:main}$
was motivated by \cite{bh1}, but the proof
presented here is completely independent.
It is algorithmic in nature,
and may be thought of as a more sophisticated version of
the Gaussian algorithm bringing an integer matrix to the Smith
normal form by elementary row and column transformations.
The extended version of Theorem~$\ref{thm:main}$ for the degenerate orbits
(corresponding to the case $n=0$)
is given in Theorem~$\ref{thm2}$.

The concept of a projective element was introduced in \cite{bh1}.
The idea is that these elements are mapped to invertible ideal classes under
Bhargava's correspondence. In the case of $\wedge^3(\zet^6)$ they were
defined via their $\SL_6(\zet)$ orbit representatives
(cf. Definition~$\ref{def:proja}$(a,b)).
It follows from our considerations that projective elements have a very convenient
description in terms of the partial derivatives of the quartic invariant
of the module (see Corollary~$\ref{cor:nabla}$). This assertion is proved for the spaces
$\mz$ as
in Theorem~$\ref{thm:main}$, but it remains valid in other spaces associated to ideal
classes in quadratic orders. (see Subsection~\ref{i:hcl}).

We also indicate a link between the Freudenthal construction and the
original Gauss's composition of quadratic forms in the Appendix at the end
of the paper.

\medskip
As a by-product of our considerations we obtain the classification of
orbits in the four ``Freudenthal" modules in the case of a field. Similar
results were obtained in \cite{clerc} and \cite[Section 5.3]{lm}
for orbits of one-dimensional subspaces in the case of complex numbers.
Our
techniques however are purely algebraic, which allows us to extend those
results to orbits of {\em elements} over arbitrary field of char$\ne 2,3$.

\medskip
\hspace{-\parindent}{\bf Theorem~$\ref{thm1}$}
{\em
\begin{enumerate}
\item[\rm (i)]
{
Let $F$ be a field of char$\ne 2,3$, let $\calg$ be the split composition algebra
$\bin, \quat, \oct$ of dimension $2,4,8$ over $F$, and let $\J=\h\calg$.
Let $(G, \m)$ be the pair (group, module) produced from $\J$
by the Freudenthal construction.

Then

\begin{itemize}
\item
There exists a $G$-invariant quartic form (the norm) on the module $\m$.

\item
The group $G$ acts
transitively on the sets of elements of rank $1,2$, and $3$ in the module
$\m$.

\item
In the case of rank\/ $4$ the group $G$ acts transitively on the set of
elements of a given norm $k$, for any $k\in F$, $k\ne 0$.
\end{itemize}
All these orbits are distinct, and the union of these orbits and $\{0\}$
is the whole module~$\m$.
}

\item[\rm (ii)]
If in addition every element of $F$ is a square,
then the same results apply to the pair $(G, \m)$ obtained
from the Jordan algebra $\J=\h F$.
\end{enumerate}
This construction yields the classification of orbits
of the irreducible representations
of simple algebraic groups listed in the following table
$$
\begin{array}{|c|c|c|}
\hline
\qquad \J \qquad       &  \mbox{Type of }G & \mbox{Highest weight of }\ \m\\
\hline
\h{\F}    &  C_3 & \omega_3\\
\h{\bin}  &  A_5 & \omega_3\\
\h{\quat} &  D_6 & \omega_5\ {\rm or}\ \omega_6\\
\h{\oct}  &  E_7 & \omega_7\\
\hline
\end{array}
$$
}

\subsection{The Freudenthal construction and higher composition laws}
\label{i:hcl}

M.~Bhargava showed that Gauss's composition law
is one in a series of at least 14 examples of the same kind (higher
composition laws) \cite{bh-phd, bhau, bh1}.
There is a certain integral linear group $G_\zet$ and
a module $V_\zet$ over it in each of his examples
such that $G_\zet$-orbits in $V_\zet$ can be described
in terms of ideal classes of orders in a number field.
Spaces underlying higher composition laws are closely related to
prehomogeneous vector spaces classified by M.~Sato and T.~Kimura \cite{SK}. In
particular, each of them is equipped with a polynomial, which is invariant
under the action of the appropriate group.

An examination of the table of higher composition laws \cite[Table 1]{bhau}
shows that each
space associated to a quadratic ring (except Gauss's composition) has a
polynomial invariant of degree four. A more detailed analysis suggests
that for each pair $(G_\zet, \ V_\zet)$ associated to a quadratic ring,
there exists a cubic Jordan algebra $\J_\zet$, such that
$(G_\zet, \ V_\zet)$ is essentially the pair produced by the Freudenthal
construction. These observations are summarized in Table 1 below.

\bigskip
\caption{The Freudenthal construction and higher composition laws}
\nopagebreak
$$
\begin{array}{|c|c|c|c|c|c|c|}
\hline
\# &
\J &
\phantom{\Bigl(}
{\dim}\ \J&
{\rm Group}\ \gr &
{\rm Rep.\ \m(\J)} &
{\dim}\ \m(\J)  &
(\lie, \gamma) \\
\hline
\hline

1 & F\phantom{\Bigl(}   & 1  & \SL_2 & {\rm Sym}^3 V_2     & 4  & G_2, \alpha_2 \\
2 & F\oplus F   & 2  & (\SL_2)^2 & V_2\tensor {\rm Sym}^2 V_2     & 6  & B_3, \alpha_2 \\
3 & \h{0}=F\oplus F\oplus F\phantom{\Bigl(}
& 3 & (\SL_2)^3& V_2\otimes V_2\otimes V_2 & 8  & D_4, \alpha_2 \\
4 & F\oplus Q_4& 5 & \SL_2\times \SL_{4} & V_2\tensor \wedge^2 V_4     & 12  & D_5, \alpha_2\\
5 & \h\bin  & 9  & \SL_6 & V(\omega _3)     & 20  & E_6, \alpha_2 \\
\hline
\hline
6 & \h\quat & 15 & D_6 & \mbox{half-spin} & 32  & E_7, \alpha_1  \\
7 & \h\oct  & 27 & E_7 & {\rm minuscule}     & 56  & E_8, \alpha_8  \\
\hline
\hline
8 & F\oplus Q_n, n\ge 3 & 1+n & {\rm SL_2\times SO}_{n+2} & V_2\tensor V(\omega_1)
& 2n+4  &{\mathfrak s\mathfrak o}_{n+6}, \alpha_2  \\
9 & \h{F}   & 6  & C_3 & V(\omega _3)     & 14  & F_4, \alpha_1 \\
\hline
\end{array}
$$

\renewcommand{\abstractname}{Summary of Table 1}
\nopagebreak
\begin{abstract}
{\bf Headers:} the table lists a semisimple
Jordan algebra $\J$ with a cubic form, its
dimension, the group $\gr$ (up to a finite subgroup or finite covering)
and its module $\mj$ produced by the Freudenthal
construction, and the dimension of $\mj$.
The last column lists a certain simple Lie
algebra $\lie$ and its simple root $\gamma$ associated to $\gr, \mj$ (see
subsection \ref{lieconn}).

{\bf Notation:} in the second column, $\h{\calg}$ stands for the Jordan
algebra of $3\times 3$-Hermitian matrices over the composition algebra
$C$; and $Q_n$ stands for the simple Jordan algebra of a quadratic form
(details are found in Subsection~$\ref{ss:ex}$). $V_n$ denotes the standard $n$-dimensional
module for $\SL_n$, and $V(\omega_i)$ denotes the simple module with
highest weight $\omega_i$ over the appropriate group.
Note that rows $2,3,4$ are special cases of row $8$ with $n=1,2,4$.

\medskip

There is a natural integral structure in each Jordan algebra $\J$ above. It
induces an integral structure in $\mj$ and $\gr$.
The first five rows of the table contain the integral group
and a module, which appear in the list of
higher composition laws
(see \cite[Table 1]{bhau} or \cite{bh1}).

Rows $6$ and $7$ do not appear in the list of higher composition laws \cite{bhau}.
However, Theorem~$\ref{thm:main}$  implies that the structure of orbits
of the projective elements is
the same as in row $\#5$. These are the first examples
of groups, more complex than a direct product of $\SL$'s,
acting in the spaces underlying higher composition laws.

Finally, the Freudenthal construction produces two more examples (rows $8$
and $9$). We conjecture that for $\#8$ the orbit structure can be
described in terms of known examples of small dimension. And the row $9$
could possibly produce a new example of a space with a composition law.
\end{abstract}

It would be interesting to describe orbits of the group $D_6(\zet)$ in the
$32$-dimensional module and of the group $E_7(\zet)$ in the
$56$-dimensional module in terms of the ideal classes of the appropriate
rings in a way similar to how this is done in \cite{bh1} for $\SL_6(\zet)$
orbits in $\wedge^3 \zet^6$.

Finally, we note that our algorithmic approach allows one to develop,
at least in some cases, a reduction theory similar to that for binary
quadratic forms.

\subsection{Connection with Lie algebras}
\label{lieconn}

Another interesting feature of the spaces underlying higher composition
laws is their remarkable connection with exceptional groups.
M. Bhargava showed that each of his pairs (group, space) may be described
in terms of the Levi decomposition of parabolic subgroups
of exceptional groups.

This connection can be stated more precisely for the Freudenthal
construction, and hence for orbit spaces associated to quadratic rings.
Here we present a more algebraic (as opposed to Lie group) realization of
this construction.
Namely, we consider a simple Lie algebra $\lie$ with a $5$-grading associated
to the minimal nilpotent orbit in $\lie$. From this setup it
is possible to extract an algebraic group $G$ and a $G$-module $M$ with
$G$-invariant quartic and skew-symmetric bilinear form, such that $(G,M)$
will produce all pairs arising from the Freudenthal construction.

The detailed description of this procedure is given below. Most of the
material of this subsection has appeared earlier elsewhere. We used
\cite[Section~2]{gw} and \cite{clerc} as the references.
For simplicity we assume in this subsection that the ground field $F$ is
the field of complex numbers $\com$.

\medskip
Let $\lie$ be a complex simple Lie algebra, and we assume that $\lie\ne
A_n, C_n$.
Let $\Phi$ be its root system with respect to a Cartan subalgebra
$\cart$; let $\Phi^+$ be a set of positive roots,
and let $\Delta$ be its collection of simple roots.
Let $\beta$ be the
highest root of $\Phi ^+$. We normalize the inner product
$\langle\cdot,\cdot\rangle$ on the real span of roots by requiring
$\langle\beta,\beta\rangle=2$. Then it follows that for any
$\alpha\in\Phi^+$, $\langle\alpha,\beta\rangle=0,1,2$; and
$\langle\alpha,\beta\rangle=2$ iff $\alpha=\beta$.

We consider the grading
\begin{equation}\label{grading}
\lie=\lie_{-2}\oplus \lie_{-1}\oplus\lie_{0}\oplus\lie_{1}\oplus\lie_{2},
\end{equation}

given by
$$
\lie_0=\cart\oplus\bigoplus_{\stackrel{\alpha\in\Phi}
{\langle\alpha,\beta\rangle=0}}
\lie_\alpha,\qquad
\lie_k=\bigoplus_{\stackrel{\alpha\in\Phi}
{\langle\alpha,\beta\rangle=k}}
\lie_\alpha\quad
k=\pm 1, \pm 2.
$$

In particular, we have $\lie_2=\lie_\beta,\
\lie_{-2}=\lie_{\beta}$.

We will also need to consider the extended Dynkin diagram of $\lie$
(with the additional vertex corresponding to the root $-\beta$).
We let $\gamma$ denote the unique simple
root whose vertex is connected to the
extended vertex $-\beta$. The standard references for root systems
and Dynkin diagrams is \cite{b1}.
The appropriate diagrams are also provided in~\cite{gw}.

The grading (\ref{grading}) has the following property:
the root space  $\lie_\alpha$ is contained in $\lie_k$
if and only if
the decomposition of $\alpha$ in the basis of simple roots $\Delta$
has coefficient $k$ at the root $\gamma$
$(k=0,\pm 1, \pm 2)$.

There exists a subalgebra $\ml\subset\lie_0$
such that $\lie_0=[\lie_{-2}, \lie_2]\oplus \ml$.
$\ml$ is a semisimple Lie
algebra, and its Dynkin diagram can be obtained from the Dynkin diagram of
$\lie$ by removing the vertex corresponding to the simple root~$\gamma$.

We let $M_\com$ be a complex connected Lie group, whose Lie algebra is
$\ml$. The group $M_\com$ and the Lie algebra $\ml$ act on the space
$\lie_1$. The resulting pairs $(M_\com, \lie_1)$ are tabulated in
\cite[Table~2.6]{gw}. An examination of that table shows that when
$\lie\ne A_n, C_n$, pairs $(M_\com, \lie_1)$ are the same (up to a finite
covering) as the pairs $(\gr, \m)$ arising from the Freudenthal construction
and listed in Table~1 above. In particular, the last column of that table
lists the simple Lie algebra $\lie$ and the simple root $\gamma$ such that
the ``Freudenthal" pair $(\gr,\m)$ is obtained from the construction
described above.

\medskip
Next, we are going show how one can describe the ``Freudenthal" quartic and
symplectic forms on $\lie_1\ (=\m)$ invariant with respect to $M_\com$.

The Lie bracket on $\lie$ induces a map $\wedge^2 \lie_1\to \lie_\beta$
which is $\lie_0$-equivariant. Since $\lie_\beta$ is one-dimensional, this
map defines a symplectic bilinear form $\{\cdot, \cdot\}$ on $\lie_1$.
In addition, we have $[\ml, \lie_2]=0$, and this implies that
$\{\cdot, \cdot\}$ is invariant with respect to $\ml$ (and hence $M_\com$).

Finally, for $X\in \lie_1$ we consider the map
$$
({\adlie\,}X)^4: \ \lie_{-2}\ \longrightarrow\ \lie_2.
$$

For fixed root vectors $x_{\pm \beta}\in \lie_{\pm 2}$ we have
\begin{equation}\label{qmap}
({\adlie\,}X)^4: \ x_{-\beta}\ \mto \ P(X)\,x_\beta \qquad
\mbox{for some }P(X)\in\com.
\end{equation}

The function $P:\lie_1\to\com$ defined by (\ref{qmap}) is in fact non-zero
homogeneous
polynomial function of degree four, invariant with respect to $\ml$ and
$M_\com$ \cite[Proposition~3.1]{clerc}. Since $\ml$ acts irreducibly on
$\lie_1$, we can normalize the polynomial $P$ so that it coincides with
the Freudenthal quartic form.

\medskip
With a little more work one can even extract the cubic Jordan algebra $\J$
from $\ml$. Namely, in each such $\ml$ there will be a simple root
$\alpha_0$, which will induce a $3$-grading on $\ml$ similar to the one we had
for $\lie$:
\begin{equation}\label{grading2}
\ml=\ml_{-1}\oplus\ml_{0}\oplus\ml_{1}.
\end{equation}

Then the space $\ml_1$ will have a structure of a Jordan algebra.

One way
to prove it is by considering a compact real form of (\ref{grading2}) and
noticing that it produces a Hermitian symmetric space, which is of tube
type (see \cite[Section 4]{clerc} for details).

The structure of a Jordan algebra on $\ml_1$
can be described in purely algebraic terms,
see \cite{tits}, or \cite[Lemma~$4$]{kac}, which is a more accessible
reference.

\medskip
To summarize, given a cubic Jordan algebra $\J$, the Freudenthal
construction produces a group $\gr$ and its module $\m$ with
$\gr$-invariant quartic and symplectic forms. Also there exists a
simple Lie algebra $\lie\ (\ne A_n, C_n)$ with grading (\ref{grading}),
which yields the same data (and $\m=\lie_1)$. Moreover the subalgebra $\ml
\subset\lie_0$ possesses grading (\ref{grading2}) such that $\ml_1$
has the structure of a Jordan algebra isomorphic to $\J$.

\subsection{Organization of the paper}
We provide basic information about split
composition algebras and their integral structures
in the beginning of Section~$2$. Then we give basic information about Jordan
algebras.
We present the Springer construction of cubic Jordan algebras and
describe the three main examples of Jordan algebras possessing
an admissible cubic form (Subsection~$\ref{ss:ex}$). In this paper
we are mainly interested in cubic Jordan algebras $\h\calg$ of
$3\times 3$ Hermitian matrices over composition algebras
(Example~$\ref{ex:jcubic}$).
In Subsection~$\ref{ss:33herm}$ we describe the
groups $\str\J$ and $\npg\J$ associated with these algebras, and describe
orbits in $\h\calg$ under the action of the norm-preserving group
$\npg{\h\calg}$ (Proposition~$\ref{prop:np-o}$).
We conclude Section~$2$ with the description of the integral structures in
Jordan algebras of Subsection~$\ref{ss:ex}$, and the description of
integral orbits in $\h\calg$ (Theorem~$\ref{np-orbits}$).
\medskip

We begin Section~3 with the description of the Freudenthal construction.
Given a cubic Jordan algebra $\J$, this construction produces a group
$\gr$ and its representation $\mj$.
We identify groups $\gr$  and its representations $\mj$
in Proposition~$\ref{prop:gener}$ and Remark~$\ref{rem:semis}$, see also
Example~$\ref{ex:sl6}$.

We present the concept of rank of elements in the module $\mj$ in
Subsection~$\ref{ss:rank}$, and we use it in the classification of
$\gr$-orbits in $\m$ in Subsection~$\ref{ss:canf}$
(Theorem~$\ref{thm1}$). The key step in our approach is the computation of
Lemma~$\ref{lcomp}$, which allows to do a complete reduction of elements of
$\m$ using purely algebraic considerations.

Most of the assertions in Subsection~$\ref{ss:canf}$ are proved under the
assumption that the ground field $F$ is an arbitrary field of
characteristic $\ne 2,3$. The restriction char\,$F\ne 2$, is essential
in the paper since it is used in many intermediate computations (see also
Remark~$\ref{badred}$). As for characteristic $3$, it appears that this
restriction may be dropped without impairing the statements. However this
assumption was made in \cite{b}, which is our main reference concerning
the Freudenthal construction, and therefore we had to incorporate it in
our considerations.

\medskip
Section~$3$ may also be viewed as a testing ground (or a simpler version) of
the techniques that we apply to studying the integral case in Section~$4$.
We introduce the integral version of the Freudenthal construction  in
Subsection~$\ref{ss:intstr}$. We develop a reduction procedure for
elements of $\mz$, which is somewhat similar to the reduction of
$n\times n$ matrices with integer entries under the elementary row and
column transformations. This is done in Lemma~$\ref{bl-fz}$, which is one
of the main technical results in Section~4. This lemma does not provide a
complete reduction of the elements in $\mz$, but it is sufficient, for
example, to describe generators of the group $\grz$ (Proposition~$\ref{genz})$
in a way similar the case of a field.

Our elementary approach does not seem to be sufficient for the complete
description of orbits in $\mz$, and we restrict our attention to the
orbits of projective elements. This concept was introduced by M.Bhargava
\cite{bh-phd, bh1}
in the context of orbit spaces associated with ideal classes in quadratic
rings. In Subsection~$\ref{ss:proj}$ we show how this concept may be
transferred to elements of the modules $\mz$. We also provide a simple
test for the projectivity of elements in terms of the quartic invariant
of the module $\mz$ (Corollary~$\ref{cor:nabla}$,
Remark~$\ref{rem:nabla}$).

In Subsection~$\ref{ss:fred}$ we explain how one can do the complete
reduction for the projective elements (Lemma~$\ref{red2}$) and summarize
all the previous results in the proof of the main result of the paper
on the structure of orbits of the projective elements
(Theorem~$\ref{thm:main}$).

In the last subsection we analyze degenerate orbits of elements of $\mz$.
We are able to obtain a complete classification of orbits of elements of
rank $1$ and $2$. These results complement results of the previous
subsection; they are stated in Theorem~$\ref{thm2}$.

\medskip
We conclude the paper with an Appendix, where we describe a link between
the Freudenthal construction and Gauss's law of composition of quadratic
forms via the Cube Law of M.~Bhargava.

\subsection{Notation and conventions}

We work over a ground field $F$ of characteristic $\ne 2,3$.
The vector space of $n\times n$ matrices over $F$ is denoted by $M_n(F)$.
The identity matrix in $M_n(F)$ is denoted by $I_n$. For an arbitrary
vector space $V$ over $F$, the symbol $\id_V$ denotes the identity linear
transformation in $\End_F(V)$.

We use the word {\em space} to denote a set with an additional structure.
Depending on a context such a space may be a vector space, a
$\zet$-module, or a set of orbits under the action of a group
(an {\em orbit space}).

Many spaces that we consider have a certain integral structure introduced
in the paper. For a space~$V$,
the corresponding integral structure will be denoted by $V_\zet$.
Very often such spaces will arise as spaces of representation
of certain algebraic groups. In such cases $G(\zet)$ will denote the
set of $\zet$-rational points of $G$ with respect to the integral
structure in $V_\zet$. This is an {\em integral form\/} of the group $G$.

We will occasionally make references to simple split finite-dimensional
Lie algebras (algebraic groups) and their root systems.
The labeling of their simple roots, fundamental weights, etc. corresponds
to that of \cite{b1}.

\bigskip
{\bf Acknowledgements.}
The author is grateful to E.~Zelmanov and M.~Racine for their help and
encouragement. Many thanks to
B. Allison,
M.~Bhargava,
O.~Loos,
E.~Neher
for helpful discussions and comments on the content of the paper.
The author is also grateful to M.~Bhargava for providing a copy of his
manuscript~\cite{bh1}.

\section{Cubic Jordan algebras and associated structures}
\label{s:jordan}

The Freudenthal construction that we will describe in Section~\ref{sec:FC}
was originally introduced in the case of the $27$-dimensional exceptional
Jordan algebra. We are going to apply this construction to other examples
of the so called cubic Jordan algebras, and we introduce the appropriate
concepts in this section.

The basic information concerning Jordan algebras is given in
Subsection~$\ref{ss:jbasic}$. However, considering the Freudenthal
construction, we will not make much use of the Jordan product $\jprod$.
Instead, we will be looking at the cubic form $N$, the trace form and
the trace bilinear form, the ``sharp" operation $\#$ and its linearization
$\times$.

All of these operations (as well as the Jordan product $\jprod$) in a cubic
Jordan algebra may be defined starting with an admissible cubic form.
We will give the appropriate construction in Subsection~$\ref{ss:sc}$,
and provide explicit examples of admissible cubic forms (and Jordan
algebras) in Subsection~$\ref{ss:ex}$.

\subsection{Composition algebras and their integral structures}
\label{ss:compint}

A (not
necessarily associative) algebra $\calg$ with a unit element
is called a {\em composition} algebra
if it has a non-degenerate quadratic form $\ncomp$ ({\em the norm})
which permits composition, i.e.,
$$%
\ncomp(xy)=\ncomp(x)\ncomp(y),\qquad
\mbox{for all }x,y\in\calg.
$$%

A famous theorem of Hurwitz states that
such an algebra is always finite-dimensional,  %
and moreover,
its dimension can only be equal to $1, 2, 4$, or $8$
(see, e.g.,  \cite[Section~II.2.6]{atoj} or \cite[Theorem~1.6.2]{SpV}).
In this paper we
restrict our attention to {\em split} composition algebras only (a
composition algebra is split, if it contains non-zero elements of zero
norm; in that case the quadratic form $\ncomp$ has maximal Witt index).
For a field $F$ there is a unique up to isomorphism split composition
algebra of a given dimension $(2,4,$ or $8)$ \cite[Theorem 1.8.1]{SpV}.
We will call it the algebra of split {\em binarions, quaternions, octonions},
and denote it by
$\bin, \quat, \oct$, respectively.
We will use symbol $\calg$ to denote an arbitrary fixed composition
algebra.

Next we will present the construction of the three split
composition algebras and show how to define an integral structure in each
of them.

Our starting point for constructing composition algebras will be the
algebra of $2\times 2$ matrices over the field $F$. This algebra has
dimension $4$ and in fact it can be taken as a model for the algebra of
split quaternions $\quat$. The quadratic form $\ncomp$ in this algebra
is the usual determinant of $2\times 2$ matrices $\det$. The composition
property for $\ncomp$ is just the multiplicativity of the determinant.

One can define the concepts of {\em trace} and {\em conjugation} of
elements in a composition algebra. In the case of quaternions they are
the trace and the symplectic involution of a $2\times 2$ matrix:
\begin{equation}
\tcomp (\bfa)=a+d,\qquad
\conj{\bfa}=\maketwo{d}{-b}{-c}a,\qquad
{\rm for}\
{\bfa }=\maketwo{a}bcd\in \quat.
\end{equation}

One can define the algebra of binarions $\bin$ to be the subalgebra of
$\quat$, which consists of the diagonal elements
\begin{equation}\label{binar}
\bin=\left\{\left.
\maketwo{a}00d\right|
a,d\in F\right\}
\end{equation}

The quadratic form, trace, and the conjugation operation in $\bin$ are
induced by those in $\quat$.

We may consider the ground field $F$ to be embedded in the binarions
\begin{equation}
F\cong \left\{\left.\maketwo{a}00a\right|a\in F\right\}
\subset\bin,
\end{equation}
and consider the induced structure of a one-dimensional composition algebra.

\medskip
Finally, we define the algebra of octonions using the Cayley-Dickson
duplication process \cite{atoj, SpV}. $\oct$ is defined to be the
direct sum $\quat\oplus\quat$, and we write an arbitrary octonion with the
aid of a formal variable $v$ (the ``imaginary unit") as
\begin{equation}
{\bfa}+\bfb v,\qquad \bfa , \bfb\in \quat.
\end{equation}

One then defines
$$
({\bf a+b}v)\cdot ({\bf c+d}v)\defeq
{\bf (ac-\conj{d}b)+(da+b\conj{c})}v,
\qquad {\bf a,b,c,d}\in \quat.
$$
\begin{equation}
\ncomp({\bf a+b}v) \defeq \det({\bf a})+\det({\bf b}),\quad
\tcomp({\bf a+b}\it v) \defeq \tcomp({\bf a}), \quad
\conj{{\bf a+b}v} \  \defeq \conj{{\bf a}}-{\bf b}v.
\end{equation}

These operations turn $\oct$ into a composition algebra.

We note that the multiplication of binarions is both associative and
commutative; the multiplication of quaternions is associative, but not
commutative; and, finally, we lose both the associative and commutative
property of the multiplication of the octonions.

We identify elements of the ground field $F$ with the appropriate
multiples of the identity element in the composition algebra. These are
the only elements which are fixed under the conjugation.

\bigskip

To define an integral structure in these algebras, we again start with
the quaternions. We define integral quaternions to be the $2\times 2$
matrices with integral entries
\begin{equation}
\quatz=
\left\{\left.
\maketwo{a}bcd\right|
a,b,c,d\in \zet\right\}.
\end{equation}

This integral structure is extended to the binarions and octonions in a
natural way
\begin{equation}
\binz=
\left\{\left.
\maketwo{a}00d\right|
a,d\in \zet\right\},\qquad
\octz=
\Bigl\{{\bfa}+\bfb v\Bigr|\ \bfa , \bfb\in \quatz
\Bigr\}.
\end{equation}

It is easy to see
that each of these $\zet$-modules has the structure of a composition
algebra (over~$\zet$), and the trace and the norm defined on them have
integer values.

\subsection{Basics of Jordan algebras}
\label{ss:jbasic}

We are going to present basic definitions and constructions concerning
Jordan algebras in this subsection. The classical reference for the
subject is N.~Jacobson's book \cite{ja}. A modern exposition of
the theory of Jordan algebras may be found in the recent monograph by
K.~McCrimmon \cite{atoj}.

\begin{defn}
{\rm
A (linear) {\em Jordan algebra} $\J$ over a field $F$ of char$\ne 2$ is a vector
space over $F$ with a bilinear product $\jprod$
(the {\em Jordan product}) satisfying the following two axioms:
\begin{equation} \label{j-ax}
x\jprod y =y\jprod x;\qquad
x^2\jprod (x\jprod y) = x\jprod (x^2 \jprod y)\qquad
{\rm \ for \ }x,y\in \J
\end{equation}
($x^2$ is defined as $x\jprod x$).
}
\end{defn}
The Jordan product is commutative by definition, but it does not have to satisfy
the associative law.

A prototypical example of a Jordan algebra is the algebra $A^+$
obtained from any associative algebra $A$ (e.g., a matrix algebra)
by defining the Jordan product to be
\begin{equation}  \label{j-assoc}
x\jprod y =\frac{1}{2}(xy+yx).
\end{equation}
Here and below the product $xy$ represents the multiplication in the
original (associative) algebra~$A$.

\medskip
We will be mostly interested in the so-called cubic Jordan
algebras, i.e., Jordan algebras,
in which every element satisfies a cubic polynomial equation.
This concept can be introduced via the concept of {\em generic norm} and
{\em generic minimal polynomial} \cite[Section VI.3]{ja}.

However in this paper we will use a different approach. Namely, we will
consider certain examples of Jordan algebras,
which are constructed from a vector
space with a cubic form. The details of the process are given in the
following subsection.

\subsection{The Springer construction of cubic Jordan algebras}
\label{ss:sc}

In this subsection we will present a construction of a class of Jordan
algebras obtained from a cubic form on a vector space, known as the
Springer construction. Our exposition will follow
\cite[Section~I.3.8]{atoj}. The original references are
\cite{spr}, \cite[Section 4]{mc-sp}. We are still working under the
assumption that char$F\ne 2,3$, and we will use it occasionally to
simplify the constructions.

\bigskip
A cubic form $\mnorm$ on a vector space $V$ over $F$ (char\,$F\ne 2,3$) is a map
$\mnorm:V\to F$ such that
\begin{itemize}
\item $\mnorm(\alpha x)=\alpha^3 \mnorm(x)$ for $\alpha\in F, x\in V$;
\item $\mnorm(x,y,z)$ is a trilinear function $V\times V \times V\to F$;
\end{itemize}
where $\mnorm(x,y,z)$ is the full linearization of $\mnorm$ given by
$$
\mnorm(x,y,z)=\frac{1}{6}
\Bigl(\mnorm(x+y+z)-\mnorm(x+y)-\mnorm(x+z)-\mnorm(y+z)+\mnorm(x)+\mnorm(y)+\mnorm(z)\Bigr).
$$
In particular, we have
$$
\mnorm (x,x,x) = \mnorm (x).
$$

We say that $c\in V$ is a {\em basepoint} for $\mnorm$, if $\mnorm(c)=1$. One then
can define the following maps
\begin{itemize}
\item
a linear map (trace) $V\to F$: $\tr(x)=3\,\mnorm(c,c,x)$;
\item
a quadratic map $V\to F$: $S(x)=3\mnorm(x,x,c)$;
\item
a bilinear map $V\times V\to F$: $S(x,y)=6\mnorm(x,y,c)$;
\item
a trace bilinear  form $V\times V\to F$: $(x,y)=\tr(x)\tr(y)-S(x,y)$.
\end{itemize}

In particular, we have
$$
\mnorm(c)=1;\qquad
S(c)=3;\qquad
\tr(c)=3.
$$

\begin{defn}
{\em
A cubic form with a basepoint $(\mnorm,c)$ on a finite-dimensional vector space
$V$ over a field $F$ of char$\ne 2,3$ is said to be {\em admissible}
or a {\em Jordan cubic}, if
\begin{enumerate}
\item[(1)]
$\mnorm$ is nondegenerate at the basepoint $c$ in the sense that
the trace bilinear form $(x,y)$ is nondegenerate;
\item[(2)]
The quadratic {\em adjoint} (or {\em sharp}) map $V\to V$,
defined uniquely by $(x^\#, y)=3\mnorm(x,x,y)$, satisfies the adjoint identity:
\begin{equation}\label{adid}
(x^\#)^\#=\mnorm(x)x.
\end{equation}
\end{enumerate}
We will also use the term {\em Jordan cubic} when referring to the
associated Jordan algebra (see Proposition~\ref{prop:jcubic}).
}
\end{defn}

The following relation holds for all $x$ in $V$:
\begin{equation}\label{trace-spur}
\tr(x^\#)=S(x).
\end{equation}

\bigskip
We define the linearization of the sharp map
\begin{equation}
x\cross y = (x+y)^\#-x^\#-y^\#.
\end{equation}
Note that
$$
x\cross x = 2 x^\#.
$$

\begin{prop}\label{prop:jcubic}
{\em \cite[Section I.3.8]{atoj}}
Every vector space with an admissible cubic form gives rise to a Jordan
algebra with unit $\one=c$ and the Jordan product given by
\begin{equation}\label{j-prod}
x\jprod y =
\frac{1}{2}\Bigl( x\cross y+\tr(x)y+\tr(y)x-S(x,y)\one \Bigr).
\end{equation}
Every element of this Jordan algebra satisfies the cubic polynomial:
\begin{equation} \label{3poly}
x^3-\tr(x)x^2+S(x)x-\mnorm(x)\one=0.
\end{equation}
We also have
$$
x^\#=x^2-\tr(x)x+S(x)\one.
$$
\end{prop}

Taking the trace of the last expression and then using~$(\ref{trace-spur})$,
one gets
$$
\tr(x^2)=\tr(x)^2-2S(x),
$$

which linearizes to
\begin{equation} \label{traces}
\tr(x\jprod y)=(x,y).
\end{equation}

\medskip
The following simple example illustrates the concepts introduced above.

\begin{example}\label{ex:mat3}
{\rm
Let $V$ be the vector space $M_{3}(F)$ of $3\times 3$ matrices over $F$.
We let the cubic form $\mnorm$ be the determinant, and we let $c$ be the
identity matrix.

Then the linear map $\tr$ is the regular trace of matrices. The sharp map
$x^\#$ produces the classical adjoint (transposed cofactor) matrix of $x$.
Equation (\ref{j-prod}) yields the Jordan product, which coincides with the
product in the Jordan algebra $M_{3}(F)^+$ given by (\ref{j-assoc}):
$$
x\jprod y=\frac{1}{2}(xy+yx).
$$
The trace bilinear form $(x,y)$ is equal to $\tr(x\jprod y)$.
It coincides in this
example with the standard trace form in the matrix algebra.
Finally, we notice that the equation (\ref{3poly}) becomes just the~
Cayley-Hamilton equation for $3\times 3$ matrices.
}
\ee
\end{example}

\subsection{Three main examples}
\label{ss:ex}

We will provide several more examples of cubic Jordan algebras in this
subsection.

The following example is a simplified version of the {\em reduced cubic
factor example} \cite[Section~I.3.9]{atoj}.
\begin{example}\label{ex:jcubic}
{\rm
Let $\calg$ be a composition algebra over a field $F$
with a quadratic form $\ncomp$ and an involution $\bar{\ }$.

We let $V$ be the space $\h{\calg}$ of $3\times 3$ Hermitian matrices over $\calg$.
An arbitrary element $A$ in $\h{\calg}$ has the form
$$
A=\makeherm{a}bc{\bf x}{\bf y}{\bf z},\qquad
\begin{array}{l}
\mbox{where }\ a,b,c \in \F, \\
\mbox{and }{\bf x,y,z}\in\calg.
\end{array}
$$
The basepoint $c$ is defined to be the identity (diagonal) matrix in
$\h\calg$, and the cubic form
$$
N:\h{\calg}\to F
$$
is given by the expression
reminiscent of the regular determinant
\begin{equation} \label{cubic1}
\mnorm(A) = \ abc-a\, {\bf x\conj{x}}-b\, {\bf y\conj{y}}-c\, {\bf z\conj{z}}+
{\bf (xy)z+\bar{\bf z}{(\bar{\bf y}\bar{\bf x})}}.
\end{equation}

The trace form $\tr$, as defined in the previous section from $\mnorm$,
coincides with the regular trace $\tr(A)=a+b+c$ and
$\btr(A,B)=\tr\Bigl(\frac{1}{2}(AB+BA)\Bigr)$.
The ``sharp" operation produces the regular adjoint matrix:
$$
A^\#=
\left(
\begin{array}{lll}

   bc-\ncomp({\bf x}) \quad                       & \conj{\bf y}\ \conj{\bf x}-c\,{\bf z} \quad   & {\bf z\, x}-b\,\conj{{\bf y}} \\
   {\bf x}\ {\bf y}-c\,\conj{\bf z}     & ac-\ncomp({\bf y})                          & {\bf \conj{z}\ \conj{y}} -a\,{\bf x}\\
   {\bf \conj{x}\, \conj{z}}-b\,{\bf y} & {\bf y \, z} -a\,{\bf x}                &  ab-\ncomp({\bf z})
\end{array}
\right).
$$

The cubic form $\mnorm$ defined by $(\ref{cubic1})$ is admissible.
}\ee
\end{example}

We note that Example \ref{ex:mat3} can be viewed as
a special case of Example \ref{ex:jcubic}
as explained in the following

\begin{rem}\label{rem:hbm3}
{\rm
If the composition algebra $\calg$ is the algebra of split binarions
$\bin$, then an arbitrary element of $\h\bin$ has the form
(see~$(\ref{binar})$)
$$
\makethree
{a}                 {\bfa_{12}}         {\bfa_{13}}
{\conj{\bfa_{12}}}  {b}                 {\bfa_{23}}
{\conj{\bfa_{13}}}  {\conj{\bfa_{23}}}  {c}
,\qquad
\mbox{where }
\bfa_{ij}=\maketwo{a_{ij}}{\scriptstyle 0}{\scriptstyle 0}{a_{ji}}\in
\bin,\quad
a,b,c, \in F.
$$

We then define a linear map $\h\bin\to M_3(F)$
\begin{equation}\label{bin3iso}
\makethree
{a}                 {\bfa_{12}}         {\bfa_{13}}
{\conj{\bfa_{12}}}  {b}                 {\bfa_{23}}
{\conj{\bfa_{13}}}  {\conj{\bfa_{23}}}  {c}
\longmapsto
\makethree
{a}       {a_{12}}  {a_{13}}
{a_{21}}  {b}       {a_{23}}
{a_{31}}  {a_{32}}  {c}.
\end{equation}

It is easy to see that this isomorphism of vector spaces
is a norm isometry, where the cubic form in $M_3(F)$ is
the regular determinant.
Hence $(\ref{bin3iso})$ defines an isomorphism of
Jordan algebras.
\ee
}
\end{rem}

\begin{example}\label{ex:quad}
{\rm
\cite[Section~$4$]{mc-sp}
Let $Q_n$ be a vector space of dimension $n$ over a field $F$,
and let $B_0$ be a
non-degenerate quadratic form on $Q_n$ with a basepoint $c_0$
(i.e., $B_0(c_0)=1$). We form a vector space $V$ by
taking the direct sum of a copy of the
ground field $F$ and $Q_n$:
$$
V=F \oplus Q_n.
$$

We then define a cubic form $\mnorm$ on $V$ by
$$
\mnorm(\alpha, x_0)=\alpha\, B_0(x_0), \qquad
\mbox{for }\alpha \in F, x_0\in Q_n.
$$
We let $c=(1,c_0)$. A simple verification
shows that $c$ is a basepoint for $\mnorm$
and that the cubic form $\mnorm$ is admissible.

In particular, the following formulas hold for $x,y \in V$:
\begin{eqnarray*}
x^\# &=& \Bigl(B_0(x_0),\ \alpha x_0^*\Bigr),\\
(x,y)&=& \alpha\beta + B_0(x_0^*, y_0),
\end{eqnarray*}
where $x=(\alpha, x_0), y=(\beta, y_0)$, $x_0^*=B_0(x_0, c_0)c_0-x_0$, and
$B_0(\cdot,\cdot)$ is a linearization of the quadratic form $B_0(\cdot)$:
$$
B_0(u,v)=B_0(u+v)-B_0(u)-B_0(v).
$$
}\ee
\end{example}

\begin{example}\label{ex:f1}
{\rm
We let $V$ be the one-dimensional vector space $V=F$. We define a cubic
form $\mnorm$ on $V$ by
$$
N(\alpha)=\alpha^3, \qquad
\alpha\in V
$$
with the obvious choice of the base point $c=1\in F$. We have
$\tr(\alpha)=3\alpha$, $(\alpha, \beta)=3\alpha\beta$,
$\alpha^\#=\alpha^2$ for $\alpha, \beta\in F$, and evidently
$N$ is an admissible cubic form.
}\ee
\end{example}

\subsection{Two groups associated to cubic Jordan algebras}
\label{ss:33herm}

Next we will introduce two groups associated to a Jordan algebra $\J$. The
definitions, due to N.~Jacobson, are valid for all finite-dimensional
Jordan (or power-associative) algebras, but we will use them only for
Jordan algebras with a cubic form.

\begin{defn}\label{def-groups}

The {\em norm-preserving group}
$$
\npg\J=\Bigl\{g\in GL(\J)\ |\  \mnorm(g A)=\mnorm(A)\mbox{ for all }A\in \J\Bigr\}
$$
is the group of all invertible $F$-linear transformations of the
vector space $\J$, which preserve the norm~$\mnorm$.

Similarly we define the {\em group of norm similarities} or
the {\em structure group}
$$
\str\J=\Bigl\{g\in GL(\J)\ |\  \mnorm(g A)=\chi(g)\, \mnorm(A)\mbox{ for all }A\in \J\Bigr\},
$$
where $\chi(g)$ is a scalar in $F$, which depends on the group element $g$
only.
\end{defn}

We have the obvious inclusion $\npg\J \subset \str\J$.

\bigskip

\begin{rem}\label{groups}
{\rm
Here we will provide a description of these groups for cubic Jordan
algebras of the form $\h\calg$
(\cite[Ch.~$14$]{sp-alg}, see also \cite{J3}, \cite[VI.7--9]{ja}).

\begin{itemize}
\item[(i)]
{\bf Case} $\J=\h{F}$.

$\str\J$ is the group of transformations of the form
$$
X\mto \gamma AXA^t,
$$

where $X\in\h{F}$, $A\in GL_3(F)$, $\gamma\in F^\mult$, $A^t$ is the transpose of $A$.

$\npg\J$ consists of transformations for which
$\gamma^3(\det A)^2=1$.

\item[(ii)]
{\bf Case} $\J=\h{\bin}$.

In this case $\J$ is isomorphic to the Jordan algebra of all $3\times 3$
matrices over $F$ (Remark~$\ref{rem:hbm3}$).

The group $\str\J$ is generated by the transformations of the form
$$
\eta(A,B): X\mto AXB^{-1}\quad
\mbox{ and }\quad
t: X\mto X^t,
$$

where $X\in M_3(F)$, $A,B\in GL_3(F)$, and $X^t$ denotes the transpose of $X$.

$\eta(A,B)$ acts trivially iff $A=B=\alpha\id, \ \alpha\in F,\alpha\ne 0$.

$\npg\J$ consists of transformations $\eta(A,B)$ for which $\det A=\det B$.

We will let $\str \J\cc$ denote the subgroup of $\str\J$ that
consists of the transformation $\eta(A,B)$. It is a subgroup of index two
in $\str\J$, and we have $\str\J\cong\str\J\cc\rtimes \{1, t\}$.
Similarly we define $\npg\J\cc=\npg\J \cap \str\J\cc$.

\item[(iii)]
{\bf Case} $\J=\h{\quat}$.

We have the isomorphisms of Jordan algebras
$\h{\quat}\cong {\cal H}(M_{6}(F), {\rm symp})$ (the $6\times 6$ symplectic
symmetric matrices) (see, e.g., \cite[p. 65]{J3}).

The group $\str\J$ is the group of transformations of the form
$$
X\mto \gamma AXA^{\rm symp}
$$

where $X\in{\cal H}(M_{6}(F), {\rm symp})$,
$A\in GL_6(F)$, $\gamma\in F$, $\gamma\ne 0$, $A^{\rm symp}$ is the transpose of $A$
with respect to the symplectic involution.

$\npg\J$ consists of transformations for which $\gamma^3 (\det A)^2=1$.

\item[(iv)]
{\bf Case} $\J=\h{\oct}$.

In this case we have that $\npg{\J}$ is a simply connected
simple algebraic group of type $E_6$, whose center is
isomorphic to the group $\bfmu _3(F)$
of third roots of unity in $F$.

The group $\str\J$ in this case is equal to the product
$\npg\J \cdot (F^\mult\id_\J)$ with $\npg\J \cap F^\mult\id_\J=\bfmu_3(F)\id$.
\end{itemize}

Later we will view groups $\str\J$ and $\npg\J$ as the group of
$F$-rational points of the algebraic groups $\astr \J$ and $\anpg\J$,
respectively. These group are connected in all cases, except $\J=\h\bin$.
For $\J=\h\bin$, the group $\astr{\J}$ (respectively, $\anpg{\J}$)
has two connected components; and the set of $F$-rational points of the
component of the identity coincides with $\str{\J}\cc$ (respectively,
$\npg{\J}\cc$).
\ee
}
\end{rem}

The following definition generalizing the concept of rank for the usual
$3\times 3$ matrices goes back to N.~Jacobson.

\begin{defn}\label{def:jrank}
Let $\J$ be a cubic Jordan algebra and let $\mnorm$ be the cubic form on $\J$.
The {\em rank} of an arbitrary element $A\in\J$ is an integer between zero and
three, which is defined by the following relations:
\begin{itemize}
\item
\rank$A =3$ iff \ $N(A)\ne 0$;
\item
\rank$A \le 2$ iff \ $N(A)=0$;
\item
\rank$A \le 1$ iff \ $A^\#=0$;
\item
\rank$A=0$ iff \ $A=0$.
\end{itemize}
\end{defn}

It is known that the rank is invariant under the action of the
groups $\npg\J$ and $\str\J$, see, e.g., \cite[Section~$2$]{J3}.

\medskip
We will conclude this subsection by describing orbits in $\h\calg$
under the action of the \np\ group $\npg{\h\calg}$.
The classification is based on the concept of rank
and is analogous to the classification of orbits in $M_3(F)$
under the action of the elementary row and column transformations.

\begin{prop}\label{prop:np-o}

\
\begin{enumerate}[\bf (i)]
\item
Let $\calg$ be the split composition algebra over a field $F$,
$\calg=\bin, \quat, \oct$, char $F\ne 2,3$.
Let $\h{\calg}$ be the cubic Jordan algebra
of $3\times 3$ Hermitian matrices over $\calg$.

Then the group $\npg{\h{\calg}}$ acts transitively on the
sets of elements of rank $1$ and $2$. In the case of rank\/ $3$,
the group $\npg{\h{\calg}}$ acts transitively on the elements
of a given norm $k$, $k\in F$, $k\ne 0$.

All these orbits are distinct and the union of these orbits and $\{0\}$ is
$\h{\calg}$.
\item
If in addition every element of $F$ is a square,
the same result holds for $\h{F}$.
\end{enumerate}
\end{prop}

It follows that the orbit representatives for the action of
$\npg{\h{\calg}}$ may be chosen to be the following
diagonal matrices:
$$
0, \quad
\makethree 1.. .0. ..0,\quad
\makethree 1.. .1. ..0,\quad
\makethree 1.. .1. ..k,\quad
k\ne 0.
$$

{\sc Proof.}

{\bf (i)}
The statement of the theorem is obvious for matrices in
$\h\bin\cong M_3(F)$ (Remark~$\ref{rem:hbm3}$).
For other split composition algebras the statement was
essentially known to N.~Jacobson, see e.g. \cite[Section~$2$]{J3}.
Alternatively, one may view this proposition as a corollary (of the proof)
of the main theorem of \cite{k}, see also Theorem~$\ref{np-orbits}$ below.

{\bf (ii)}
The proof of (i) does not work in the case of Hermitian matrices over the
ground field $F$. The result here depends on the arithmetic properties
of the ground field. When $F$ satisfies the assumptions of (ii), the
assertion is well known, see, e.g. \cite[Theorems~$6.5, 6.6$]{jba}.
\ep

\subsection{Operations and identities in cubic Jordan algebras}
\label{ss:id}

Given an arbitrary Jordan algebra $\J$, one can define the
{\em Jordan triple product}
\begin{equation}\label{jtp}
\{x,y,z\}  =  (x\jprod y)\jprod z+ x\jprod (y\jprod z)-(x\jprod z)\jprod
y,\qquad
x,y,z\in\J.
\end{equation}

When $\J\subseteq A^+$ for some associative algebra $A$, the Jordan triple
product has a simple expression in terms of the associative operation in
$A$:
\begin{equation}\label{a:triple}
\{x,y,z\}=\frac{1}{2}(xyz+zyx),\qquad
x,y,z\in \J\subseteq A^+.
\end{equation}
This operation is very important for generalizations of Jordan algebras.
We mentioned it here, since we need to introduce another operation
$V_{x,y}:\J\to\J$ defined by
\begin{eqnarray}\label{a:v}
V_{x,y} (z) & = & \{x, y, z\},\qquad
x,y,z\in \J.
\end{eqnarray}
We will use the operation $V_{x,y}$ in
Subsections~$\ref{ss:rank}$ and~$\ref{ss:if}$ in the definition
of rank of elements of the module~$\mj$.

\medskip

Next we list several identities which relate the triple product, the cross
product, and the trace bilinear form in an arbitrary {\em cubic} Jordan
algebra~$\J$.

It follows from $(\ref{3poly})$ that
$$
X\jprod X^\#=\mnorm(X)\,\one
$$

for any element $X$ in $\J$.

One can linearize this identity
(see, e.g., \cite[Section~${2}$]{J3}, \cite[Section~5.2]{SpV})
and get

\begin{eqnarray} \label{a:id1}
\{X,Y,X\}+2Y\cross (X^\#)&=&(X,Y)X\\ \label{a:id2}
(X,Z)Y+(Y,Z)X&=&2\{X,Z,Y\}+(X\cross Y)\cross Z\\ \label{a:id3}
N(Y)X+(X,Y^\#)Y&=&(X\cross Y)\cross Y^\#\\ \label{a:id4}
(X,Z^\#)X+\{Z,X^\#,Z\}&=&(X\cross Z)^\#.
\end{eqnarray}

One can define the * operation
for an arbitrary element $s$ of the \np\ group
$\npg\J$ by the relation:
\begin{equation} \label{a:adj}
\Bigl(\ s   (X),\ Y\ \Bigr) \ =
\ \Bigl(\ X,\ s^* (Y)\ \Bigr) \qquad \mbox{for any }X,Y\in \J.
\end{equation}

It satisfies the following identities:
\begin{equation} \label{a:crossid}
{s^{*-1}}(X\cross Y)= s (X) \cross s (Y) \qquad \mbox{for any }s\in\npg\J,\ X,Y\in\J.
\end{equation}
\begin{equation} \label{a:tripleid}
s \Bigl(\{X,Y,Z\}\Bigr)= \Bigl\{\ s (X),\ {s^{*-1}}(Y),\ s (Z)\ \Bigr\} \qquad
\mbox{for any }s\in\npg\J,\ X,Y,Z\in\J.
\end{equation}

\subsection{The integral case and orbits under the \np\ group}
\label{ss:33int}

The integral structure $C_\zet$
of the composition algebras $C$ in Subsection~$\ref{ss:compint}$ induces an
integral structure in the spaces $\h{C_\zet}$.

We note that $\h{C_\zet}$ is not closed under the Jordan product
because of the factor $\frac{1}{2}$, but
it is easy to see from the expressions given in Example~$\ref{ex:jcubic}$,
that the cubic form $\mnorm$, the trace $\tr$ and the trace bilinear form
$\btr(\cdot, \cdot)$ take values in (and onto) $\zet$. We also have that
$\h{C_\zet}$ is closed under the sharp operation $\#$, and
$X\cross Y\in\h{C_\zet}$ for any $X,Y\in \h{C_\zet}$.

We consider the \np\ groups defined in Subsection~$\ref{ss:33herm}$,
and we look at
their subgroups of elements which preserve the integral submodule
$\h{C_\zet}$. This subgroup is an integral form of the appropriate group.
We consider the action of each of these groups on the space of
integral Hermitian matrices.

The structure of the orbits in $\h{C_\zet}$ under the action
of the \np\ groups is described in
Theorem~\ref{np-orbits} below. Before stating the theorem
we will recall the following well-known

\begin{defn}
{\rm
We say that an $n\times n$ matrix $A$ is in the {\em Smith normal form},
if $A$ is a diagonal matrix
$$
A=\diag{d_1, d_2}{\dots}{d_n},
$$
where
$$
\mbox{all }d_i\mbox{'s are integers, }\quad
d_i \,|\, d_{i+1},\quad d_i\ge 0\quad
\mbox{ for }1\le i \le n-1
$$
and all zeros on the diagonal are located in the lower right corner.
}
\end{defn}

\begin{theorem} \label{np-orbits}
Let $\h{\calg_\zet}$ be the $\zet$-module of $3\times 3$ Hermitian matrices
over the split composition {\em ring} $\calg_\zet$, $\calg_\zet=\bin_\zet,
\quat_\zet, \oct_\zet$.

Then every element of $\h{\calg_\zet}$ is equivalent to an element in the
{\em Smith normal form} under the action of the group $\npg{\h{\calg_\zet}}$.
Distinct elements in the Smith normal form lie in distinct
$\npg{\h{\calg_\zet}}$-orbits.
\end{theorem}

{\sc Proof.}
The assertion of this theorem when $C_\zet=\oct_\zet$
was proved in \cite{k} (we used the term {\em canonical diagonal form}
to represent Smith normal form).

The reasoning of that paper also applies in the case of integer
quaternions and was stated as Corollary there \cite[p.~294]{k}.
This assertion may also be stated in terms of integer skew-symmetric
matrices (cf. Remark~\ref{groups}(iii)), see \cite[Theorem~IV.1]{newman}.

In case of integer binarions, one could also repeat the argument of
\cite{k} to arrive to the same conclusion. An alternative way
to prove the theorem is to notice that the isomorphism of
Remark~$\ref{rem:hbm3}$ yields a $\zet$-linear norm isometry of
$\zet$-modules
$\h{\binz}$ and $M_{3}(\zet)$
A consequence of this fact is that the action of the \np\
group in $\h\binz$ can be expressed in terms of the elementary row and columns
transformations of the regular $3\times 3$ matrices over $\zet$.
Hence in the case of
integral binarions the assertion of the theorem is equivalent to the
assertion on the Smith normal form for the usual $3\times 3$ matrices over
$\zet$.
\ep

\begin{rem}
{\rm
An important feature of the theorem is that we have here a chain of
embedded spaces
\begin{equation}
\h{\bin} \subset \h{\quat} \subset \h{\oct},
\end{equation}
and the action of the corresponding groups there.
Theorem~\ref{np-orbits} gives us
a uniform description of orbits for all three spaces
in terms of diagonal matrices contained in each of these spaces.
\ee
}
\end{rem}

We also have the following trivial corollary to Theorem~\ref{np-orbits}
(which also applies in the case of general $n\times n$ matrices over
$\zet$ and their orbits under the elementary row and column transformations).

\begin{cor}
Let $\h{\calg_\zet}$ and $G_\zet$ be as in Theorem~$\ref{np-orbits}$,
and let $n$ be an integer, $n\ne 0$.
The group $G_\zet$ acts transitively on set of matrices of determinant $n$
if and only if $n$ is a squarefree integer.
\end{cor}

\section{The Freudenthal construction}
\label{sec:FC}

\subsection{Preliminary facts on the module $\mj$ and the group $\gr$}
\label{ss:fc}

In this section we consider a certain class of modules and
linear groups acting on them.
Historically, the first example of this kind was introduced by
H.~Freudenthal
in \cite{fr14} in the process of constructing the $56$-dimensional
representation of the group $E_7$ from the $27$-dimensional exceptional
Jordan algebra.
These modules were studied axiomatically (under the name of the Freudenthal
triple systems) in \cite{b}, \cite{fau}, \cite{fer}.

We will consider examples of Freudenthal triple systems
of the form $\mj$, where $\J$ is a cubic Jordan algebra
of Section~$\ref{s:jordan}$.
We will say that the module $\mj$ and its automorphism group $\gr$ (see below)
are obtained from $\J$ using the~{\bf Freudenthal construction}.

\medskip
We will use \cite{b} as the main reference for this subsection.
We consider a vector space $\m=\mj$ constructed from the space $\J$ in the
following way
\begin{equation}\label{eq:mj}
\mj=F\oplus F \oplus \J \oplus \J,
\qquad \mbox{ where $\J$ is a cubic Jordan algebra over $F$.}
\end{equation}

We have $\dim \m=2\, \dim\J+2$, and an arbitrary element $x$ of the space
$\m$ has the form
\begin{equation} \label{gen}
x=\mel\alpha \beta A B,\qquad
{\rm where\ } \alpha, \beta \in F,\quad A,B \in \J.
\end{equation}

We have the following skew-symmetric bilinear and quartic forms
on $\m$ defined by
$$
\{ x,y\}= \alpha\delta - \beta\gamma+ (A,D)-(B,C)
$$
$$
q(x)= 8 (A^\#, B^\#) -8\alpha \mnorm(A) - 8\beta \mnorm(B)
-2\Bigl((A,B)-\alpha\beta\Bigr)^2.
$$
Here we have
$x=\mel\alpha \beta A B$, $y=\mel\gamma \delta C D$.
$(\cdot, \cdot)$ is the trace bilinear form,
$\mnorm$ is the norm, and $\#$ is the
sharp map in the Jordan algebra $\J$
(see Subsections~$\ref{ss:sc}$,~$\ref{ss:ex}$ for detail).

Later on, when we turn to the integral case, it will be more convenient
for us to consider the modified form $q'$:
\begin{equation}\label{qmod}
q'(x)= -4 (A^\#, B^\#) +4\alpha \mnorm(A) +4\beta \mnorm(B)
+\Bigl((A,B)-\alpha\beta\Bigr)^2,
\qquad q=-2q'.
\end{equation}

We will often refer to the form $q$ as the {\em norm form} or just the
{\em norm} in the module~$\m$.

We can linearize the form $q$ and get a symmetric four-linear form
$q(x,y,z,w)$ such that $q(x,x,x,x) = q(x)$. It follows that both the bilinear
and four-linear form are non-degenerate. And hence we can define a
trilinear operator $T:\m \times \m \times \m \to \m$ by the following
rule: for given $x,y,z \in \m$, $T(x,y,z)$ is the unique element in $\m$
such that
\begin{equation} \label{deft}
\{ T(x,y,z), w\} = q(x,y,z,w)\qquad \mbox{for any $w\in \m$}.
\end{equation}

\begin{defn}
The group $\gr$ is defined to be the group of all invertible $F$-linear
transformations of $\m$ that preserve these forms, i.e.,
\begin{equation} \label{invar}
\Bigl\{ \sigma (x),\sigma (y)\Bigr\} = \{ x,y\},\qquad q\Bigl(\sigma (x)\Bigr) =
q(x)\quad
\end{equation}
for any $\sigma \in \gr$.
\end{defn}

It follows immediately from $(\ref{deft})$ and $(\ref{invar})$ that for
any $\sigma \in \gr$
\begin{equation} \label{autot}
T\Bigl(\sigma (x), \sigma (y), \sigma (z)\Bigr)=
\sigma \Bigl( T(x,y,z) \Bigr).
\end{equation}

We have the following four types of transformations in the group $\gr$:
\begin{itemize}

\item[]
For any $C\in \J$

\begin{equation} \label{phi}
\phi(C)\quad : \quad \mel \alpha \beta A B\  \mto
\mel {\alpha+(B,C)+(A, C^\#)+\beta \mnorm(C)} {\quad    \beta}
 {\quad A+\beta C} {\quad  B+A\cross C +\beta C^\#}.
\end{equation}

\item[]
For any $D\in \J$

\begin{equation} \label{psi}
\psi(D) \quad : \quad
\mel \alpha \beta A B\  \mto
\mel \alpha {\quad \beta + (A,D) + (B, D^\#) + \alpha \mnorm(D)}
{A+ B\cross D+\alpha D^\#} {B+\alpha D}.
\end{equation}

\item[]
In addition, for every norm similarity $s\in \str\J$
(cf. Definition~$\ref{def-groups}$)
we have the transformation $\tf s \in \gr$ defined by:
\begin{eqnarray} \label{t}
\tf s& : &
\mel \alpha\beta A B\  \mto
\mel {\lambda ^{-1}\alpha}{\quad \lambda\beta}
     {\quad s(A)}{\quad s^{*^{-1}}(B)}.
\end{eqnarray}
Here $\lambda \in \F$ is such that $\mnorm(s(x))=\lambda \mnorm(x)$,
and $s^*$ is the linear transformation adjoint to $s$ with respect
to the trace bilinear form $(\cdot, \cdot)$.
We will mostly
use such transformations when $s\in \npg\J$ and so $\lambda =1$.

\item[]
Finally, we have the transformation $\tau$, which acts by:
\begin{equation} \label{tau}
\tau :  \quad
\mel\alpha\beta A B \ \mto \ \mel {-\beta} {\alpha} {-B} A.
\end{equation}
\end{itemize}

The transformations $\phi(\cdot)$ and $\psi(\cdot)$ are conjugate to each
other by $\tau$:
$$
\tau \phi(C) \tau^{-1}=\psi(-C).
$$
In addition we have the following relations in $\gr$:
\begin{equation}\label{comm-rel}
\begin{array}{rcl}
\phi(-\one)\psi(\one)\phi(-\one) &=& \tau\\
\tau^2 &=& -\id_\m\\
\tau T(s) &=& T(s^{*-1})\tau\\
T(s)\phi(C) &=& \phi(\lambda^{-1}\, s(C))\, T(s)\\
T(s)\psi(C) &=& \psi(\lambda\, s^{*-1}(C))\, T(s).
\end{array}
\end{equation}
Here $C\in\J$ and $s\in \str\J$ satisfies
$\mnorm(s(D))=\lambda \mnorm (D)$.

Note that when matrices $C$ and $D$ above have rank 1, the transformations
$\phi$ and $\psi$ have the following simpler form:
\begin{eqnarray}\label{t1s}
\phi(C) & : &
\mel\alpha\beta A B\  \mto\
\mel{\alpha+(B,C)}{\quad    \beta}{\quad A+\beta C}{\quad  B+A\cross C}\\
     \label{t2s}
\psi(D) & : &
\mel \alpha \beta A B\  \mto \
\mel \alpha{\quad \beta + (A,D)}{\quad A+ B\cross D}{\quad    B+\alpha D}.
\end{eqnarray}

\medskip
The following proposition gives a more precise description of the group
$\gr$ and its representation $\mj$. We will let $\astr\J$ (respectively,
$\ag\m$) denote the algebraic group whose group of $F'$-rational points is
$\astr{\J\tensor_F F'}$ (respectively, $\ag{\m\tensor_F F'}$) for all extension
fields $F'$ of $F$. The symbol $\aJ$ will denote the vector (algebraic)
group defined by $\J$, i.e., $\aJ(F')=\J\tensor_F F'$ taken with the
additive group structure.

\begin{prop}\label{prop:gener}
\
Let $\J$ be a Jordan algebra $\h{\calg}$ with $\calg=F,\bin, \quat, \oct$
over a field $F$ (char$F=p\ne 2,3)$
and let $\m=\mj$.

\begin{itemize}

\item[{\rm (i)}]

Then the group $\gr$
is generated by elements $\phi (C)$, $\psi (D)$, $\tf s$.

\item[{\rm (ii)}]
The  group $\gr$ is the set of $F$-points
of an absolutely almost simple linear algebraic group $\ag\m$,
which is defined over $F$ and $F$-split.

The group $\gr$ is connected (except the case $\J=\h\bin$). It has
a two-element center, and its quotient modulo the
center is a simple group of adjoint type.

In the case $\J=\h\bin$ the same result is true for the connected component
$\ag\m\cc$, which is a subgroup of index $2$ in $\gr$.

\item[{\rm (iii)}]
The following table lists the types of the group $\ag\m$ as well
as the highest weight of its irreducible representation on the space $\m$.

$$
\begin{array}{|c|c|c|}
\hline
\qquad \J \qquad       &  \mbox{Type of }\ag\m & \mbox{H.w. of }\ \m\\
\hline
\h{\F}    &  C_3 & \omega_3\\
\h{\bin}  &  A_5 & \omega_3\\
\h{\quat} &  D_6 & \omega_5\ {\rm or}\ \omega_6\\
\h{\oct}  &  E_7 & \omega_7\\
\hline
\end{array}
$$

\end{itemize}

\end{prop}

{\sc Proof.}

(i) was proved in \cite[Theorem~3]{b}. The theorem is stated there when $\J$
is $27$-dimensional exceptional Jordan algebra, but it remains true for
other algebras in the list. We will give a somewhat different proof
of this statement when we describe the generators of the group $\gr$ in
the integral case (Proposition~$\ref{genz}$).

\medskip
(ii)
It is not hard to see that algebraic equations defining the group $\ag\m$
have integer coefficients.
Then reducing modulo $p$, we can assume that it is defined over the prime
field $F_p$, and hence over~$F$.

The $F$-split torus in $\ag\m$ is the image of the diagonal split torus of
$\astr\J$ under the mapping
$T:\astr\J\to\ag\m$, see (\ref{t}). Such a torus has the
``right rank", and it was described explicitly in Remark~\ref{groups}
(the case of $\J=\h\oct$ was considered in \cite[Theorem 3.5]{gar}).

Next we notice that ${\phi(\aJ)}, {\psi(\aJ)}, T{(\astr\J)}$
are closed subgroups
of $\ag{\m}$, since they are homomorphic images of the algebraic
groups $\aJ, \aJ, \astr\J$, respectively.
These groups are connected when
$\J=\F, \quat, \oct$ (see Remark~\ref{groups}).
Since by (i) the group $\ag\m$ is generated by these subgroups,
it is connected in these cases.

If $\J=\bin$, the group ${\astr\J}$ has two connected components, and in fact
is isomorphic to the semi-direct product ${\astr\J}\cc\rtimes C_2$, where
$C_2$ is the group of order $2$ generated by the element $t$
corresponding to the transpose operation (Remark~\ref{groups}(ii)).
This fact and commutation relations~$(\ref{comm-rel})$ imply that
$\ag\m$ may have one or two components. The {\em algebraic} subgroup $\alg{H}$,
generated by ${\phi(\aJ)}, {\psi(\aJ)}, T({\astr\J}\cc)$,
is connected and closed. In addition, it is known
(see, e.g., \cite[Section~$7.5$]{hum})
that $\alg H$ is
generated by ${\phi(\aJ)}, {\psi(\aJ)}, T({\str\aJ}\cc)$ as an
abstract group.
The analysis of the action of the group $\gr$
in the $20$-dimensional module $\mj$ (cf.
Example~\ref{ex:sl6}) shows that the element $T(t)$ does not lie
in $\alg H$.
Hence we have a decomposition into two cosets
$\ag\m=\alg H \cup T(t)\alg H$, which implies that $\alg H$ is the connected
component of $\ag \m$.

The group $\gr$ contains the central element $\tau^2$, which acts as
the scalar $-1$
on the module $\m$. The statement about the center and the simplicity
of the quotient modulo the center was proved in \cite[Theorem~6]{b}
(the theorem was stated there for
$\J=\h\oct$, but it remains true for other cases as well).

\medskip
(iii) If char$F=0$,
the group $\ag\m$ may be identified by its Lie algebra, which is
the Tits-Kantor-Koecher construction of $\J$ \cite{koe}.

An analysis of a slightly different version of the Freudenthal construction
may be found in \cite[$2.22-2.26$]{sp-alg}
with the resulting groups being identified in Section
$14.31$ of the same book.

The corresponding irreducible representation can be identified by its
highest weight vector and its dimension. The highest weight vector in such
a representation may be chosen to be $(1,0,0,0)$ in $\mj$.
\ep

\begin{rem} \label{rem:semis}
{\rm
Another example of an admissible cubic Jordan algebra
$$
\J=F\oplus Q_n, n\ge 1
$$
was given in Example~$\ref{ex:quad}$.

In this case the Freudenthal construction
(assuming the quadratic form has maximal Witt index)
produces the semi-simple $F$-split
algebraic groups of type $\SL_2\times \SO_{n+2}$ acting on the
tensor product $V_2 \tensor V(\omega_1)$ of the
irreducible $\SL_2$-module $V_2$ of dimension two and the irreducible
$\SO(n+2)$-module of highest weight $\omega_1$ (and dimension $n+2$).

For small $n$ the action of the group $\gr$ in $\m$ is essentially isomorphic
to the action of the group $G$ in $\m$ listed in the table below.

$$
\begin{array}{|l|l|l|l|}
\hline
n & G & \m & {\rm Comment}\\
\hline
1 & \SL_2\times \SL_2 & V_2\tensor {\rm Sym}^2 V_2 &
using\ Spin_3\cong \SL_2\\
\hline
2 & (\SL_2)^3 & V_2\tensor V_2\tensor V_2 &
using\ Spin_4\cong \SL_2\times \SL_2\\
\hline
4 & \SL_2\times \SL_4 & V_2\tensor  \wedge^2(F^4) &
using\ Spin_6\cong \SL_4\\
\hline
\end{array}
$$

These cases are listed separately in Table $1$.
}\ee
\end{rem}

\begin{rem}\label{badred}
{\rm
We note that when char\,$F=2$, the quartic form $q'$ will reduce to
$$
q'(x)=(\alpha\beta-(A,B))^2,\qquad
x\in \m.
$$
Hence the group $\gr$ will become the group of transformations preserving
the symplectic form $\{\cdot, \cdot\}$ and the quadratic form
$q'_2=(\alpha\beta-(A,B))$. Next we notice that the linearization of
$q'_2$ will produce exactly the symplectic form $\{\cdot, \cdot\}$.
Hence the group $\gr$ in the case char\,$F=2$ coincides with the orthogonal
group on the vector space $\m$ with respect to the quadratic form $q'_2$.
}\ee
\end{rem}

We will conclude this subsection with an explicit description of the
module $\m=\mj$ and the group $\gr$ for the Jordan algebra
$\J=M_3(F)\cong \h\bin$.
We will let $\gr\cc$ denote the subgroup generated by
$\phi(\J), \psi(\J), T(\str\J\cc)$ for $\J=M_3(F)$
(cf. Remark~$\ref{groups}$(ii)). This is a subgroup of index two in $\gr$,
and it coincides with the group of $F$-rational points of the connected
component of the algebraic group $\ag\m$ (cf. Proposition~$\ref{prop:gener}$).

\begin{example} \label{ex:sl6}
{\rm
We let $\J=M_3(F)$, and
we are going to describe how the action of the group $\gr$
in the $20$-dimensional module $\m=\mj$ is related to the action
of the group $\SL_6(F)$ in the space $\wedge^3 F^6$.

Note that the group $\SL_6(F)$ has a center isomorphic to $\bfmu_6(F)$,
and the quotient $\SL_6(F)/\bfmu_3(F)I_6$ acts faithfully in
$\wedge^3 F^6$. We are going to explicitly describe an isomorphism of
vector spaces $\theta: \mj\to \wedge^3 F^6$ and an isomorphism of groups
$\theta':\gr\cc\to \SL_6(F)/(\bfmu_3(F)I_6)$ satisfying:
\begin{equation}\label{ex:sl3iso}
\theta(g\cdot v)=\theta' (g)\cdot \theta(v)\quad
\mbox{ for }g\in \gr\cc, v\in\mj.
\end{equation}
With a slight abuse of notation we will still use $6\times 6$ matrices to
represent elements of the quotient $\SL_6(F)/\bfmu_3(F)I_6$.

We begin by introducing the following notation. Let
$$
\{e_1, e_2, e_3, f_1, f_2, f_3\}
$$
be a standard basis of $F^6=F^3\oplus F^3$. Next we introduce the
following ``dual" set of linearly independent vectors in $\wedge^2 F^6$:
$$
e_1^*=e_2\wedge e_3,\quad
e_2^*=e_3\wedge e_1,\quad
e_3^*=e_1\wedge e_2,
$$
$$
f_1^*=f_2\wedge f_3,\quad
f_2^*=f_3\wedge f_1,\quad
f_3^*=f_1\wedge f_2.
$$

Then we define a correspondence $\theta$ between bases of $\m$ and
$\wedge^3F^6$ by
\begin{eqnarray}\label{ex:sl1}
\theta: \hspace{2cm}
\mel 1000 & \mapsto & e_1\wedge e_2 \wedge e_3,\nonumber \\
\mel 0100 & \mapsto & f_1\wedge f_2 \wedge f_3,\\
\mel 00{E_{ij}}0 & \mapsto & e_i\wedge f_j^*,\nonumber\\
\mel 000{E_{ij}} & \mapsto & f_i\wedge e_j^*,
\hspace{4cm} 1\le i,j\le 3.
\nonumber
\end{eqnarray}

Next we define the homomorphism $\theta'$. First we will do it assuming
that every element of $F$ is a cube of another element of $F$.
The group $\gr\cc$ is generated by transformations $\phi(\J), \psi(\J),
T(\str\J\cc)$.
First we define the map $\theta'$ for $\phi(\J), \psi(\J)$ by
\begin{equation}\label{ex:sl2}
\theta':\quad
\phi(A)\mto
\maketwoblock   {I_3}{0}A{I_3},\qquad
\psi(B)\mto
\maketwoblock   {I_3}B{0}{I_3},
\end{equation}
where $A,B\in\J=M_3(F)$, and each block in
$\maketwoblock\cdot\cdot\cdot\cdot$ represents
a matrix in $M_3(F)$.

Next we define a homomorphism
$$
\theta'': GL_3(F)\times GL_3(F)\to \SL_6(F)/\bfmu_3(F)I_6
$$
by
\begin{equation}\label{ex:sl3}
\theta'':
(A,B)\ \mapsto\
\zeta_A\zeta_B
\maketwoblock{A/\det (A)}0
0{B/\det (B)},\qquad
A,B\in \GL_3(F), \zeta_A, \zeta_B\in F^\mult,
\end{equation}
where $\zeta_A, \zeta_B$ are chosen so that
$$
\zeta_A^3=\det(A), \quad
\zeta_B^3=\det(B).
$$

Since the map $\theta''$ has values in $\SL_6(F)/\bfmu_3(F)I_6$,
the relation $(\ref{ex:sl3})$ does not depend on the choice of the third
roots $\zeta_A, \zeta_B$. We also note that $\ker\theta''=F(I_3, I_3)$.
Recall that $\str\J\cc \cong GL_3(F)\times GL_3(F)/F^\mult(I_3, I_3)$
(cf. Remark~$\ref{groups}$(ii)),
and hence the map $\theta''$ produces a well defined homomorphism
\begin{equation}\label{ex:sl31}
\theta': T(\str\J\cc) \to \SL_6(F)/\bfmu_3(F)I_6.
\end{equation}

A direct computation shows that the generators of $\gr\cc$ and
$\SL_6(F)/\bfmu_3(F)I_6$
associated via $(\ref{ex:sl2}),(\ref{ex:sl3}), (\ref{ex:sl31})$
define the
same linear transformations in the (isomorphic) vector spaces
$\m$ and $\wedge^3 \F^6$.
Hence the maps $\theta, \theta'$ define an ``isomorphism" of the pairs
$(\gr, \m)$ and $(\SL_6(F)/\bfmu_3(F)I_6,\ \wedge^3 F^6)$ in the sense
of~$(\ref{ex:sl3iso})$.

In the case of an arbitrary field $F$, the quantities $\zeta_A,
\zeta_B$ are elements of a suitable field extension of $F$.
However, since $\zeta_A^3, \zeta_B^3$ are elements of $F$ and
the matrix in~$(\ref{ex:sl3})$ acts in $\wedge^3 F^6$ (not just in $F^6$),
the expression~$(\ref{ex:sl3})$ will produce a well defined $F$-linear
transformation of $\wedge^3 F^6$.
\ee
}
\end{example}

\subsection{The rank of elements of $\m$}
\label{ss:rank}

\begin{defn}\label{def:rank}
Let $\m=\m(\J)$, where $\J$ is a cubic Jordan algebra. The {\em
rank} of an element $x\in \m$ is an integer between $0$ and $4$,
which is uniquely defined by the following relations:
\begin{itemize}
\item
rank $x = 4$  iff \ $q(x)\ne 0$;
\item
rank $x \le 3$  iff \ $q(x)=0$;
\item
rank $x \le 2$  iff \ $T(x,x,x)=0$;
\item
rank $x \le 1$  iff \  $3\,T(x,x,y)+\{x,y\}\,x=0$
for all $y\in \m$;
\item
rank $x = 0$ iff \ $x=0$.
\end{itemize}
\end{defn}

The expressions defining rank $2,3,$ and $4$ are quite natural; they
appeared elsewhere before, see, e.g., \cite{clerc}. %
However our coordinate-free
relation defining $\rank x\le 1$ appears to be new in this
context.

The following lemma is an immediate corollary of the definitions of $q(x)$
and $T(x,x,x)$.

\begin{lemma}\label{lem:rank}
The rank of elements is preserved under the action of the group $\gr$.
\end{lemma}

The expression defining the rank of $x$ are given by homogeneous
algebraic equations in terms of the coordinates of $x$. We list
the appropriate polynomials in the following

\begin{rem}\label{rem:rp}
{\rm
These polynomials will be described in terms of their value at an
arbitrary element $x\in \m=\mj$ of the form
\begin{equation}\label{elem1}
x=\mel\alpha \beta{\ A}B.
\end{equation}

\begin{itemize}
\item
The {\em quartic rank polynomial} is the quartic form $q$:
\begin{equation} \label{p4}
q(x)= 8 (A^\#, B^\#) -8\alpha \mnorm(A) - 8\beta \mnorm(B)
-2[(A,B)-\alpha\beta]^2.
\end{equation}

\item
It was computed in \cite{b} that
\begin{eqnarray}  \label{p3} \label{txxx}
T(x,x,x)        \nonumber
&=& \biggl( -\alpha ^2 \beta + \alpha\, (A,B)
- 2 \mnorm(B), \qquad
\alpha \beta ^2 -\beta\,(A,B) + 2 \mnorm(A),  \\
&& \phantom{-}
2B \cross A^\#-2\beta\, B^\# -[(A, B) - \alpha\beta]A,\\
&&                 \nonumber -2A\cross B^\#+2\alpha\, A^\#
+[(A,B)-\alpha\beta]B \biggr).
\end{eqnarray}

\item
We define {\em quadratic rank polynomials} to be the following expressions:
\begin{eqnarray} \label{p2}
\alpha A-B^\#, \quad \beta B - A^\# ,\qquad Q(x),
\end{eqnarray}
where $Q$ is a quadratic polynomial function with values in End$(\J)$,
i.e., $Q(x)$ is a linear transformation of $\J$ defined by
\begin{equation}\label{q}
Q(x) : \quad C \quad \mto\quad (\alpha\beta - (A,B)) C +2\, V_{A,B}(C),
\qquad \quad C\in \J.
\end{equation}
where $V_{A,B}$ was defined in $(\ref{a:v})$.
\ee
\end{itemize}
}
\end{rem}

\begin{lemma}\label{lem:rank1}
Let $x=\mel\alpha\beta A B$ be an element of $\m=\m(\J)$.
The quadratic rank polynomials
\begin{equation} \label{e1r1}
\alpha A-B^\#, \quad \beta B - A^\#,\quad Q(\bx) \qquad
\mbox{ are equal to zero}
\end{equation}
if and only if
\begin{equation} \label{e2r1}
3\,T(x,x,y)+\{x,y\}\,x \quad
\mbox{ is equal to zero \ \ for all }y\in \m.
\end{equation}
\end{lemma}
{\sc Proof.}
It was computed in \cite[p.~88]{b} that linearization of $T(x,x,x)$ yields
\begin{eqnarray*}  \label{p3lin}
3\,T(x,x,y)
&=&
\Bigl( -2\alpha \beta \gamma -\alpha^2 \delta + \gamma\, (A,B)
+\alpha\, (C,B) +\alpha\, (A,D) - 6 \mnorm(B,B,D), \qquad\\
&& \quad\
\gamma \beta ^2 +2\alpha\beta\delta-\delta\,(A,B)-\beta\, (C,B)
-\beta\, (A,D)+ 6 \mnorm(A,A,C),  \\[2mm]
&&\ \
-[(C, B)+(A,D)-\alpha\delta - \beta\gamma]A
-(A,B)C +\alpha\beta C-\\
&&\ \ %
- 2\delta B^\# -2\beta B\cross D
+2 D \cross A^\#+2B\cross(A \cross C),\\[2mm]
&&\quad\
[(C,B)+(A,D)-\alpha\delta - \beta\gamma]B+
(A,B)D -\alpha\beta D+\\
&&\quad\
+2\gamma A^\#+2\alpha A\cross C
-2C\cross B^\#-2A\cross(B\cross D)
 \Bigr)
\end{eqnarray*}
for $x=\mel\alpha \beta A B,\ y=\mel \gamma\delta C D$.

We can rewrite it using the identity ~$(\ref{a:id3})$
and the definition of the form $\{\cdot, \cdot\}$ as
\begin{equation} \label{r1eq}
3\,T(x,x,y)+\{x,y\}\,x =
\end{equation}
\vspace{-8mm}
\begin{eqnarray*} \label{r1int-main}
&=&
\Bigl(
-[3\alpha\beta-(A,B)]\gamma+2(\alpha A -B^\#, D),\\
&& \quad\
[3\alpha\beta-(A,B)]\delta-2(\beta B -A^\#, C),\\
&& \quad\
[3\alpha\beta-(A,B)]C- 2 (\beta B-A^\#)\cross D
+2(\alpha A-B^\#)\delta-2Q(x)(C),\\
&& \ \,
-[3\alpha\beta-(A,B)]D+ 2 (\alpha A-B^\#)\cross C
-2(\beta B-A^\#)\gamma+2Q(x')(D)
\Bigr).
\end{eqnarray*}

Here $x'=\mel\beta \alpha B A$.

\medskip
It follow from the last expression (using that char\,$F\ne 2$) that
$$
3\,T(x,x,y)+\{x,y\}\,x=0 \quad
\mbox{for any }y\in \m
$$

if and only if
$$
\alpha A-B^\#=0, \quad \beta B - A^\# =0,\quad
Q(x)=0,\ Q(x')=0, \quad
3\alpha\beta-(A,B)=0.
$$

This proves the ``$\Leftarrow$" implication of the lemma.

\medskip

Now suppose (\ref{e1r1}) holds. To prove ``$\Rightarrow$",
we need to show that $Q(x')=0$ and $3\alpha\beta-(A,B)=0$.

We start with the latter one. We have $Q(x)(C)=0$ for any $C\in\J$.
In particular, this is true for $C=\one$.
We then compute using~(\ref{q}) and~$(\ref{traces})$:
$$
0=\tr(Q(x)(\one))=3\alpha\beta-3(A,B)+2\tr(A\jprod  B)=
3\alpha\beta-3(A,B)+2(A, B)=3\alpha\beta-(A,B).
$$

To prove $Q(x')=0$ we use the definition of the Jordan triple
product~$(\ref{jtp})$ to notice that
$$
Q(x)(C)-Q(x')(C)=2Q(x)(\one)\jprod C,
$$
and the statement again follows from the fact that $Q(x)(C)=0$
for any $C\in \J$.
\ep

\subsection{The canonical form in the case of a field}
\label{ss:canf}

\begin{lemma} \label{bl-f}
Let $\J$ be as in Proposition~$\ref{prop:gener}$ and $\m=\mj$.
Every non-zero element of the module $\m$ can be brought to the form
$$
\mel 1 \beta {\ A} 0
\qquad for\ some\ \beta\in \F,\ A\in\J
$$
by an appropriate element of $\gr$.
\end{lemma}

{\sc Proof.} We start with an arbitrary element $x_1$ of the form
$\mel {\alpha_1}{\beta_1} {A_1} {B_1}$.

First we show that we can transform it to an element in which the
component $B_1$
is not zero. If $B_1\ne 0$, there is nothing to do. If $A_1 \ne
0$, the transformation $\tau$ does the trick.
Now let us assume $A_1= B_1 = 0$ in $x_1$. After application of $\tau$ if
necessary, we may assume that $\alpha _1\ne 0$. Then applying $\psi (D)$
with any non-zero $D$, we get an element $x_2$ in which the
component at the position $B_1$ is not zero.

Thus we proceed with an element $x_2$ of the form $(\alpha _2,\beta_2,
A_2, B_2)$ with $B _2\ne 0$. If $\alpha _2\ne 1$ we argue as follows.
Since elements of $\J$ of rank $1$ span the whole $\J$, and the trace form
in $\J$ is non-degenerate, we may find an element $C\in\J$ of rank 1, such
that $(B_2,C)\ne 0$. Scaling $C$, we may assume that $\alpha_2+(B_2,C)=1$.
The condition ${\rm rank}\,C=1$ implies $C^\# =0$ and $N(C)=0$. After
applying the transformation $\phi (C)$ to $x_2$, we get an element $x_3$
in which the first component is equal to $\alpha _2+(B_2,C)=1$.

We arrived at the element $x_3$ of the form
$$
\mel 1 {\beta_3}  {A_3} {B_3}.
$$

The application of $\psi(D)$ with $D=-B_3$ brings this element to the
desired form. \ep

\begin{lemma}
\label{lcomp} An element
$$
\mel \alpha  \beta {\ \diag{a_1}{a_2}{a_3}} 0, \quad\alpha\ne 0
$$
can be mapped by a transformation in $\gr$ to the elements
\begin{itemize}
\item[\rm(i)]
\hspace{2cm} $
\mel \alpha {\ \beta-2\frac{a_2a_3}{\alpha}\,c}
{\quad \diag{a_1+\beta c-\frac{a_2a_3}{\alpha}\, c^2}{\ a_2}{\ a_3}} 0,
$
\item[\rm(ii)]
\hspace{2cm} $
\mel \alpha {\ \beta-2\frac{a_1a_3}{\alpha}\,c}
{\quad \diag{a_1}{\ a_2+\beta c-\frac{a_1a_3}{\alpha}\, c^2}{\ a_3}} 0, $
\item[\rm(iii)]
\hspace{2cm} $
\mel \alpha {\ \beta-2\frac{a_1a_2}{\alpha}\,c}
{\quad \diag{a_1}{\ a_2}{\ a_3+\beta c-\frac{a_1a_2}{\alpha}\, c^2}} 0, $
\end{itemize}
where $c$ is an arbitrary element in $\F$.

This lemma is also valid in the case $F=\zet$ assuming
$\alpha|a_2, \alpha|a_3$ in {\rm (i)};
$\alpha|a_1, \alpha|a_3$ in {\rm (ii)};
$\alpha|a_1, \alpha|a_2$ in {\rm (iii)}.
\end{lemma}

{\sc Proof.} The proof of this lemma is a direct computation. The action
of
$$
\phi (C),\ C=\diag{c}00
$$
and then
$$
\psi(D),\ D=-\frac{1}{\alpha}\, \diag0{a_3c}{a_2c}
$$
gives the first element. The appropriate modifications of these yield the
other two elements. \ep

\begin{rem} \label{rem1}
{\rm
When the component $B$ of an element of the form $(\ref{elem1})$ in $\m$ is
equal to zero, the relations $(\ref{p4})-(\ref{p2})$ in the definition of
the rank have the following simpler form
\begin{itemize}
\item
rank $x \le 1$  iff
\begin{equation} \label{r1s}
\alpha A=0,\quad A^\# =0,\quad \alpha\beta=0;
\end{equation}
\vspace{-7mm}
\item
rank $x \le 2$  iff
\begin{equation} \label{r2s}
\alpha ^2\beta=0, \quad \alpha\beta ^2 +2N(A)=0, \quad \alpha\beta
A=0,\quad \alpha A^\# =0;
\end{equation}
\vspace{-7mm}
\item
rank $x \le 3$  iff
\begin{equation} \label{r3s}
8\alpha N(A)+2\alpha ^2\beta ^2=0\quad (i.e., q(x)=0);
\end{equation}
\vspace{-7mm}
\item
rank $x = 4$  iff
\begin{equation} \label{r4s}
q(x)\ne 0.
\end{equation}
\end{itemize}
}\ee
\end{rem}

\begin{theorem} \label{thm1}
\

\begin{enumerate}
\item[\rm (i)]
{
Let $F$ be a field of char\ $\ne 2,3$, let $\calg$ be the split composition algebra
$\bin, \quat, \oct$ over $F$, and let $\J=\h\calg$.
Let $(G, \m)$ be the pair (group, module) produced from $\J$
by the Freudenthal construction.

Then

\begin{itemize}
\item
There exists a $G$-invariant quartic form (the norm) on the module $\m$.

\item
The group $G$ acts
transitively on the sets of elements of rank $1,2$, and $3$ in the module
$\m$.

\item
In the case of rank\/ $4$ the group $G$ acts transitively on the set of
elements of a given norm $k$, for any $k\in F$, $k\ne 0$.
\end{itemize}
All these orbits are distinct, and the union of these orbits and $\{0\}$
is the whole module~$\m$.
}

\item[\rm (ii)]
If in addition every element of $F$ is a square,
then the same results apply to the pair $(G, \m)$ obtained
from the Jordan algebra $\J=\h F$.
\end{enumerate}
This construction yields the simple algebraic groups and their irreducible
representations described in Proposition~$\ref{prop:gener}$.
\end{theorem}

{\sc Proof.}

The $G$-invariant quartic form $q$ defined on the space $\m$ is a part of
the Freudenthal construction. One defines the rank of elements of $\m$ using the
quartic form (Definition~$\ref{def:rank}$).

We are going to show that the group $\gr$ acts transitively on the set of
elements of a given rank/given norm. Namely, we will show that every
element of $\m$ is $\gr$-equivalent to one of the following elements
in the ``canonical" form:
\begin{eqnarray}
\label{r1}
\rank 1 & : & \mel 10 {\diag 000} 0,\\
\label{r2}
\rank 2 & : & \mel 10 {\diag 100} 0,\\
\label{r3}
\rank 3 & : & \mel 10 {\diag 110} 0,\\
\label{r4}
\rank 4 & : & \mel 10 {\diag 11k} 0, \quad k\in F, k\ne 0.
\end{eqnarray}
These elements have distinct rank/norm,
and therefore by Lemma~$\ref{lem:rank}$ and the definition of $\gr$
they lie in distinct orbits.

We start with an arbitrary non-zero element $x$ in $\m$. By
Lemma~\ref{bl-f} we can transform $x$ to the element
\begin{equation} \label{semican}
x_1=\mel 1 \beta {\ A} 0.
\end{equation}
We have rank\,$x$=rank\,$x_1$, and we can use Remark~\ref{rem1} when
computing the rank of $x_1$. We also make use of the rank of elements
of $\J$ (Definition~$\ref{def:jrank}$).
\medskip

If rank\,$x_1=1$, then it follows from Remark~\ref{rem1} that $\beta=0$
and $A=0$. So the element $x_1$ is already in the form~$(\ref{r1})$.

If rank\,$x_1=2$, then it follows from Remark~\ref{rem1} that $\beta=0$.
The condition $\alpha A^\# =0$ of~$(\ref{r2s})$ implies rank\,$A=1$.
By Proposition~$\ref{prop:np-o}$
there exists $s\in \npg\J$ such that
$$
s(A)=\diag{1}00.
$$
The action of the transformation $\tf{s}\in \gr$ on the element $x_1$
brings it to the form~$(\ref{r2})$ as desired.
\medskip

Let us now consider the case of rank\,$x_1=3$. We start by showing that we
can assume that the $\beta$-component of $x_1$ is zero. If $\beta\ne 0$,
it follows from $(\ref{r3s})$ that
$$
N(A)=-\frac{1}{4}\beta ^2\ne 0.
$$
Again, by Proposition~$\ref{prop:np-o}$
we can find an element $s\in \npg\J$ that brings $A$ to the diagonal
form:
$$
s(A)=\diag{1}1n\quad {\rm with\ } n=N(A).
$$
And hence we have
$$
\tf{s}(x_1)=\mel 1 \beta {\ \diag{1}1n} 0.
$$
Taking $c=\beta/2$ in Lemma~\ref{lcomp}(iii) we see that we can map
the above element to the element
$$
x_2=\mel 1 0 {\ A'} 0
$$
with component $\beta$ equal to zero.

We still have $\rank x_2=3$, and
it follows from $(\ref{r3s})$ that $N(A')=0$. On the other hand, $A^\#
\ne 0$, since the opposite would imply that rank\,$x_2 \le 2$. Hence
rank\,$A'=2$. Again we can find an $s'\in \npg\J$ that brings $A'$ to the
diagonal form:
$$
s'(A')=\diag{1}10,
$$
and hence
$$
\tf{s'}(A')=\mel 1 0 {\ \diag{1}10} 0.
$$
This is an element of the form~$(\ref{r3})$,
and we are done with the case of rank 3.

\medskip
The last case is the case of rank $4$.

We are still working with an element $x_1$ of the form $(\ref{semican})$.
Since $\rank x_1=4$, we have $q(x_1)\ne 0$.

There are three possible subcases for the element
$$
x_1=\mel 1 \beta {\ A} 0.
$$

\medskip
{\em Subcase $1$.} $\beta =0$.

In this subcase the expression for the norm becomes
$$
q(x_1)=-8\, N( A).
$$
Since $q(x_1)\ne 0$, we have $N(A)\ne 0$, and we can find an $s\in \npg\J$
that brings $A$ to the form $\diag{1}1k$, $k\ne 0$. It follows that
$\tf{s}(x_1)$ has the form~$(\ref{r4})$ as desired.

\medskip
{\em Subcase $2$.} $\beta \ne 0$ and rank\,$A\ge 2$.

We start by bringing $A$ to the diagonal form of Proposition~$\ref{prop:np-o}$
by an appropriate
element $s\in \npg\J$. Then $s(A)$ has the form $\diag{1}1n$, $n\in \F$. We
can apply Lemma~\ref{lcomp}(iii) with $c=\beta /2$ to the element
$\tf{s}(x_1)$, and get a new element with $\beta$-component equal to zero.
This brings us to the subcase 1, and we are done.

\medskip
{\em Subcase $3$.} $\beta \ne 0$ and rank\,$A\le 1$.

We again start by bringing $A$ to the diagonal form by an
appropriate element $s\in \npg\J$. Then $s(A)$ has the form
$\diag{\epsilon}00$, $\epsilon=0$ or $1$ depending on the rank\,$A$. We
can apply Lemma~\ref{lcomp}(ii) and (i) (if necessary) with $c=1/\beta$,
and obtain element
$$
x_2=\mel \alpha \beta {\diag 110} 0.
$$ This brings us to
the subcase 2, and we are done again.

We showed that an arbitrary non-zero element can be brought to one of the
elements in the canonical form~$(\ref{r1})-(\ref{r4})$.
Hence the union of the orbits described in the hypothesis
is the whole~$\m$.
\ep

\section{The Integral version of the Freudenthal construction}
\subsection{The integral structure in the module $\mj$}
\label{ss:intstr}

In this section we consider the integral structure $\mz$ in the module
$\m$ and the integral form $\grz$ of the group $\gr$ that acts on it.
We define
\begin{equation}\label{mz}
\mz\defeq\m(\jz)=\zet\oplus\zet\oplus\jz\oplus\jz,
\end{equation}
where $\jz$ is one of
\begin{equation} \label{zlist}
\jz:\qquad
\h{\octz},\  \h{\quatz},\  \h{\binz},\
\h{\zet},\  \zet\oplus \zet \oplus \zet,
\end{equation}
see Subsection~$\ref{ss:33int}$

Most of the definitions and assertions of this section apply also in the
cases
$$
\zet\ \mbox{ and }\ \zet\oplus Q_n(\zet) \quad
\mbox{(cf.
Remark~\ref{rem:semis})}.
$$
However these cases require somewhat different treatment and
we will not consider them in this paper.

An arbitrary element of the module $\mz$ has the form
$$
\mel\alpha\beta AB,
\qquad
\mbox{where } \alpha, \beta \in \zet \mbox{ and } A,B \in \jz.
$$

Since $\jz$ is closed under $\#$, and $\tr$ and $\btr(\cdot, \cdot)$ are
integer-valued,
we have that the quartic form $q$, the modified form $q'$
(see \ref{qmod})), and the skew-symmetric form~$\{\cdot,\cdot\}$ have integer
values on $\mz$. The expression~$(\ref{txxx})$ implies
that $T(x,x,x)\in \mz$ for $x\in\mz$.

\begin{defn}
{\rm
We define the {group $\grz$} to be the group of invertible
$\zet$-linear transformations of $\mz$ that preserve the quartic
form $q$ and the bilinear form $\{\cdot, \cdot \}$.
}
\end{defn}
It is immediate that transformations
\begin{equation} \label{elz}
\phi(C), \psi(D),\ \tf{s},\ \tau \quad
{\rm with}\ C,D\in \jz,\ s\in \str{\jz}
\end{equation}
preserve $\mz$, and therefore lie in the group $\grz$.

\begin{defn}
{\rm
We will call the transformations of the form $(\ref{elz})$ the
{\em elementary\/} transformations of $\mz$.
}
\end{defn}

If $F$ is the field of characteristic zero, then
the group $\grz$ may also be defined as the set of integral points
of the group $\gr$ with respect to the integral structure on the
module $\m$ determined by (\ref{mz}).

We will prove in Proposition~$\ref{genz}$ below that $\grz$ is generated by
the elementary transformations~$(\ref{elz})$.
We will make use of this assertion in the next examples, which give an
explicit description of $\mz$ and $\grz$ in two special cases.

\begin{example} \label{ex:sl6int}
{\rm
This example gives a  description of the integral version
of Example~$\ref{ex:sl6}$.
Here we have
$\jz=M_3(\zet)\cong\h\binz$ and
$\mz=\m(\jz)=\zet\oplus\zet\oplus\jz\oplus\jz$.
By analogy with Example~$\ref{ex:sl6}$, we define the group $\grz\cc$ to
be the subgroup of $\grz$ generated by
$\phi(\jz)$, $\psi(\jz)$, $T(\str\jz\cc)$.
We are going to establish the isomorphism of the integral pairs
$(\grz\cc, \jz)\cong(\SL_6(\zet), \wedge^3\zet^6)$ in the sense
of~$(\ref{ex:sl3iso})$

We keep the notation for bases and define the isomorphism $\theta$ of
$\zet$-modules $\mz$ and $\wedge^3 \zet^6$ as in Example~$\ref{ex:sl6}$.

We have $\bfmu_3(\zet)=\{1\}$ and hence
$\SL_6(\zet)=\SL_6(\zet)/\bfmu_3(\zet)I_6$ acts faithfully in
$\wedge^3 \zet^6$.
We are going to define the map
$\theta':\grz\cc\to\SL_6(\zet)$ in a way similar to Example~$\ref{ex:sl6}$.

First we define $\theta'$ for $\phi(\jz), \psi(\jz)$ by
\begin{equation}\label{ex:sl2int}
\theta':\quad
\phi(A)\mto
\maketwoblock   {I_3}{0}A{I_3},\qquad
\psi(B)\mto
\maketwoblock   {I_3}B{0}{I_3},
\end{equation}
where $A,B\in\jz=M_3(\zet)$.

As for the group $\str\jz\cc$, we notice that
it consists of the transformations of the form
$$
X\mapsto AXB^{-1}
$$
with $X\in\jz, A,B\in \GL_3(\zet)$ and $(A,B)$ acting trivially iff
$A=B=\pm I_3$ (cf. Remark~$\ref{groups}$(ii)).

Next we will define a homomorphism
$$
\theta'': \GL_3(\zet)\times \GL_3(\zet)\to \SL_6(\zet).
$$
Since in this case the determinants of the matrices involved are $\pm 1$,
we can define $\theta''$ by the following simpler formula:
\begin{equation}\label{ex:sl3int}
\theta'':
(A,B)\ \mapsto\ \maketwoblock{(\det B)A}0
0{(\det A)B},\qquad
A,B\in \GL_3(\zet).
\end{equation}

The map $\theta''$ factors through the quotient
$\GL_3(\zet)\times \GL_3(\zet)/\{\pm(I_3, I_3)\}$ and yields a well-defined
map
\begin{equation}\label{ex:sl31int}
\theta': T(\str\jz\cc)\to \SL_6(\zet).
\end{equation}

A direct computation shows that the generators of the groups $\grz\cc$ and
$\SL_6(\zet)$
associated via $(\ref{ex:sl2int}),(\ref{ex:sl3int}), (\ref{ex:sl31int})$
define the
same linear transformations in the (isomorphic) $\zet$-modules
$\mz$ and $\wedge^3 \zet^6$.
Hence the maps $\theta, \theta'$ define an ``isomorphism" of the pairs
$(\grz\cc, \mz)$ and $(\SL_6(\zet),\ \wedge^3 \zet^6)$ as desired.

It will be proved in Proposition~$\ref{genz}$ that every element in $\grz$
is a product of the elementary transformations~$(\ref{elz})$. Commutation
relations~$(\ref{comm-rel})$ imply that every element in $\grz$ has
the form $g_0$ or $T(t) g_0$, where $g_0\in \grz\cc$ and $t$ is the
transpose operation, see Remark~$\ref{groups}$(ii).
Hence we have $\grz=\grz\cc \cup T(t)\grz\cc$.
\ee
}
\end{example}

\begin{example}\label{ex:FFF}
{\rm
Let $\J=F\oplus F \oplus F$ and hence $\jz=\zet\oplus\zet\oplus\zet$,
$\mz=\m(\jz)=\zet\oplus \zet\oplus \jz\oplus \jz$.
In this example we are going to explicitly describe the relation between
the action of the group $\grz$ in the $8$-dimensional module $\mz$ and the
action of the group $\SL_2(\zet)\times \SL_2(\zet)\times \SL_2(\zet)$ in the space
$\zet^2\tensor \zet^2 \tensor \zet^2$.

We will use the more compact notation $\Gamma$ to denote the
group $\SL_2(\zet)\times \SL_2(\zet)\times \SL_2(\zet)$.
Note that the group
$\Gamma$
has a center that consists of the matrices $\{(\pm I_2, \pm I_2, \pm I_2)\}$,
and the quotient
$\Gamma/K_4$ acts faithfully in
$\zet^2\tensor \zet^2 \tensor \zet^2$, where
$$
K_4=\{(I_2, I_2, I_2),(I_2, -I_2, -I_2),(-I_2, I_2, -I_2),(-I_2, -I_2, I_2),\}.
$$

Let $G_0$ be the subgroup in $\grz$ generated by transformations
$\phi(\jz)$, $\psi(\jz)$, and $T(\pm \id_\jz)$.
We are going to describe explicitly an isomorphism of $\zet$-modules
$\theta:\mz \to \zet^2\tensor\zet^2\tensor\zet^2$ and an isomorphism
of groups
$\theta':G_0 \to\Gamma/K_4$ satisfying
\begin{equation}\label{ex:F0}
\theta(g\cdot v)=\theta' (g)\cdot \theta(v)\quad
\mbox{ for }g\in G_0, v\in\mz.
\end{equation}
With a slight abuse of notation we will use triples of $2\times 2$
matrices to
represent elements of the quotient $\Gamma/K_4$.

First we will establish an isomorphism of $\zet$-modules
$\mz\cong\zet^2\tensor\zet^2\tensor\zet^2$.

Let $\{E_1, E_2, E_3\}$
be the standard basis of $\jz=\zet\oplus\zet\oplus\zet$,
and let $\{e_1, e_2\}$ be the standard basis of $\zet^2$.
We define $\theta$ on bases and extend it by $\zet$-linearity
\begin{eqnarray} \label{ex:F1}
&&
\mel 1000 \mapsto e_1\tensor e_1\tensor e_1,\quad\ \
\mel 0100 \mapsto e_2\tensor e_2\tensor e_2,\nonumber \\
&&
\mel 00{E_1}0 \mapsto e_1\tensor e_2\tensor e_2,\quad
\mel 00{E_2}0 \mapsto e_2\tensor e_1\tensor e_2,\quad
\mel 00{E_3}0 \mapsto e_2\tensor e_2\tensor e_1,\quad\\
&&
\mel 000{E_1} \mapsto e_2\tensor e_1\tensor e_1,\quad
\mel 000{E_2} \mapsto e_1\tensor e_2\tensor e_1,\quad
\mel 000{E_3} \mapsto e_1\tensor e_1\tensor e_2.\nonumber
\end{eqnarray}

Next we define the map $\theta'$ on generators of $G_0$
(and onto generators of $\Gamma/K_4$) by
\begin{eqnarray} \label{ex:F2}
\phi(a_1, a_2, a_3) & \mapsto&
\Biggl(
{%
\maketwo{1}{a_1}{0}{1},
\maketwo{1}{a_2}{0}{1},
\maketwo{1}{a_3}{0}{1}}
\Biggr)\nonumber \\
\psi(a_1, a_2, a_3) & \mapsto &
\Biggl(
\maketwo{1}{0}{a_1}{1},\maketwo{1}{0}{a_2}{1},\maketwo{1}{0}{a_3}{1}
\Biggr)\\
T(\pm \id_\jz) & \mapsto &
\pm
\Bigl(
I_2, I_2, I_2
\Bigr),\nonumber
\end{eqnarray}
where $(a_1, a_2, a_3)\in \jz$.

A direct computation shows that the generators of $G_0$ and
$\Gamma/K_4$
associated via $(\ref{ex:F2})$
define the
same linear transformations in the (isomorphic) $\zet$-modules
$\mz$ and $\zet^2\tensor \zet^2\tensor \zet^2$.
Hence the maps $\theta, \theta'$ define an ``isomorphism" of the pairs
$(G_0, \mz)$ and $(\Gamma/K_4,\ \zet^2\tensor \zet^2\tensor \zet^2)$
in the sense of~$(\ref{ex:F0})$.
}\ee
\end{example}

\subsection{Reduction in the module $\m(\jz)$}
\label{ss:red}

\begin{defn}
{\rm
Given a free $\zet$-module $M$,
we say that an integer $d$ {\em divides} an
element $x\in M$ if $x=dx'$ for some $x'\in M$.

We define the $\gcd$ of a collection of elements in
$M$ to be the greatest integer that
divides all these elements. By definition the $\gcd$ is  a positive
integer, if at least one of the elements is nonzero.
}
\end{defn}

\begin{defn}\label{prim}
{\rm
An element $x$ of a free $\zet$-module $M$ is said to be {\em primitive},
if $\gcd (x)=1$.
}
\end{defn}
Evidently, for any non-zero element $x$ we have
$
x = \gcd(x)\ x',
$
where $x'$ is primitive.

\begin{lemma} \label{gcd}
Let $M$ be a free $\zet$-module and
let $g$ be an element in {\rm End}${}_\zet(M)$ such that $g^{-1}$ exists and
also lies in {\rm End}${}_\zet(M)$.
Then for any $x$ in $M$ and any non-zero integer\/ $d$
$$
d\ divides\ x\qquad
\mbox{\rm if and only if}\qquad
d\ divides\ g(x).
$$
\end{lemma}

{\sc Proof.}
Obvious.

\bigskip

\begin{defn}\label{def:reduced}
{\rm
An element $x$ of the module $\mz=\m (\jz)$ is said to be {\em reduced},
if it is of the form
\begin{equation} \label{reduced}
x=\mel\alpha\beta{\ A}0,\qquad
\alpha,\beta\in \zet,\ A\in\jz,
\end{equation}
with $\alpha>0,\quad \alpha |\beta, \ \alpha |A$. We say that $x$ is a
{\em diagonal reduced} element, if in addition $A$ is a diagonal matrix.
}
\end{defn}

We note that for a reduced element $x$ as in (\ref{reduced}) we have
$\gcd(x)=\alpha$. And a reduced $x$ is primitive if and only if
$\alpha=1$.

\begin{lemma}[\sc Reduction Lemma I]  \label{bl-fz}
\

\hspace{-\parindent}%
Let $\jz$ be one of $\h{\binz}, \h{\quatz}, \h{\octz}, \zet\oplus\zet\oplus\zet$.
Every non-zero element ${\bx}$ of the module $\mz=\m (\jz)$
is equivalent to a {\em diagonal reduced} element under the action of
a series of elementary transformations
\begin{equation} \label{elem-tr}
\phi (C), \psi(D), \tf{s}, \tau \quad
with \ C, D \in \jz,\ s\in \npg{\jz}.
\end{equation}
\end{lemma}

{\sc Proof.}

{\bf (a)}
First we are going to prove the assertion when $\jz$
is one of $\h{\binz}, \h{\quatz}, \h{\octz}$.

We are going to apply transformations $\phi (C), \psi(D), \tf{s}$
to obtain a reduction
procedure similar to the Gaussian elimination process for
regular integer matrices.

The proof consists of two parts.
We describe four reduction steps in the first part, and in the second part
we show how we use these steps to bring an arbitrary non-zero element to a
reduced form.

{\bf Part I.}

We start the proof of the lemma by explaining how we can use the
elementary transformations~(\ref{elem-tr}) to get the reduction.
We will use the
relations $(\ref{phi})$, $(\ref{psi})$, and the ability to bring any
element in $\jz$ to the (diagonal) Smith normal form by an element
in $\npg \jz$ (see Theorem~\ref{np-orbits}).
Note that transformations $s\in \npg\jz$ preserve the norm, and hence
we have $\lambda=1$ in (\ref{t}).

The four reduction steps are summarized below:
\medskip

(RED1): \
$\mel\alpha\beta{\ A}B \longrightarrow \mel{g_1}**0$
with $0<g_1\le\ \min \{|\alpha|,\gcd( B)\}$
\medskip

(RED2): \
Assume $\alpha>0$, \ $\alpha \not | \gcd(\beta, A)$;

\ \phantom{(RED2): }%
$\mel\alpha\beta{\ A}0 \longrightarrow \mel{*}{\beta'}{A'}{{}*}$
with $0<\min\{\beta', \gcd(A')\}<\alpha$

\bigskip
(RED3): \
$\mel\alpha\beta{\ A}B \longrightarrow \mel *{g_2}0{{}*}$
with $0<g_2\le \min \{|\beta|,\gcd(A)\}$

\medskip
(RED4): \
Assume $\beta>0$, \ $\beta \not | \gcd(\alpha, B)$;

\ \phantom{(RED4): }%
$\mel\alpha \beta{\ 0} B \longrightarrow \mel{\alpha'} *{*} {B'}$
with $0<\min\{\alpha', \gcd(B')\}<\beta$

\medskip
The symbol ``$*$" above represents some element at the appropriate
position. We will not be interested in that entry for the moment.

The idea of (RED1) is to use the interaction of the elements $(\alpha, B)$
in the left column to replace them by $(g_1, 0)$. (RED2) is designed to be
applied after (RED1), and its idea is to (eventually) reduce the
pair $(\beta, A)$ in the right column modulo $g_1$.

(RED3) and (RED4) are mirror images of (RED1) and (RED2) with the roles of
the left and right columns reversed.

We now proceed to describing these steps in detail.
\medskip

(RED1): \
$\mel\alpha\beta{\ A} B \longrightarrow \mel{g_1} * * 0$
with $0<g_1\le\ \min \{|\alpha|,\gcd( B)\}$

We start with an element $x\in \mz$ of the form
\begin{equation}
x=\mel\alpha\beta{\ A} B,\qquad
\alpha,\beta\in \zet,\ A,B\in\jz.
\end{equation}

We can use the transformation $\psi (D)$ with an appropriately chosen $D \in
\jz$ to reduce entries of the element $B$ modulo $\alpha$. This
essentially means that the $9, 15,$ or $27$ integers in the entry $B+\alpha D$ of
$\psi(D)({\bx})$ will lie in the interval $0 \dots |\alpha | -1$.

Conversely, we can reduce the entry $\alpha$ modulo $m_2$, where
$m_2 = \gcd (B)$.
Without loss of generality we can assume that $B$ is in the Smith normal
form (to get it, first find an $s\in \npg\jz$ such that $s(B)$ is in the
Smith canonical form; then replace $x$ by $\tf{s^{*^{-1}}}\ ({\bx}))$;
the value of $\alpha$ will not be affected by this). This remark
implies that the component $B$ of the (new) element $x$ will have the form
$$
B=\diag{m_2}**\qquad{\rm with}\ m_2=\gcd(B).
$$
Next, we find
an integer $c_1$ such that $\alpha +m_2 c_1$ is between $1$ and $m_2$,
and let $C$ be the diagonal matrix $\diag{c_1}00$.
The ``$1,1$"-entry of the element $\phi (C) ({\bx})$ is
$$
\alpha+(B,C)+(A, C^\# ) +\beta \mnorm(C)
$$
and since rank $C$=1, it is equal to
$$
\alpha+(B,C)=\alpha+m_2 c_1,
$$
which is in $\{1,\dots, m_2\}$,
and we have reduced $\alpha$ modulo $\gcd(B)$
(note that we chose the set of residues,
which excludes $0$, but contains $m_2$, since we want to keep this
component non-zero).

We can now reduce the entries $\alpha$ and $B$ modulo each other, and
continue to do so as long as the entry at the position $B$ is non-zero.
The value of the component $\alpha$ will decrease after each pair of
iterations, while remaining positive. This condition guarantees that the
process will terminate at some step. This procedure
will bring $x$ to the desired form, and since at each step we reduced
modulo $\alpha$ or $\gcd(B)$, the condition on $g_1$ is satisfied. This
completes the description of (RED1).

\bigskip

(RED2): \
Assume $\alpha>0$, \ $\alpha \not | \gcd(\beta, A)$;

\ \phantom{(RED2): }%
$\mel \alpha\beta A0 \longrightarrow \mel *{\beta'}{A'}{{}*}$
with $0<\min\{\beta', \gcd(A')\}<\alpha$.

We start with an element $x$ as stated in the assumption.
We can apply a transformation $\tf{s}$ to bring the entry $A$ to the
Smith normal form without changing $\alpha$ and $\beta$, so without
loss of generality we assume that $A$ is already in the Smith normal form
$A=\diag{a_1}{a_2}{a_3}$, $\ a_1|a_2,\ a_2|a_3,\ \alpha>0$.

Let us now take matrix $D=\diag{0}1{d_3}$ (we will choose the value for
$d_3$ later). Then, since $B=0$ and $\mnorm(D)=0$, we have
$$
\psi(D)({\bx})=
\mel\alpha {\ \beta +(A,D)}{\quad  A+\alpha D^\# }{\ \alpha D}.
$$

We have
$D^\# =\diag{d_3}00$, and hence the ``$1,2$"-entry of the above element
is
$$
\diag{a_1+d_3 \alpha}{\quad a_2}{\quad a_3}.
$$
We can find $d_3\in \zet$ that brings $a'_1=a_1+d_3 \alpha$ into the
integer interval $1,\dots,\alpha$.
If $a'_1<\alpha$, then we are done. And if $a'_1=\alpha$, it means
$\alpha | A$. Then
we can use the element $a'_1$, ideas from (RED3), and the non-divisibility
assumption (which translates now into $\alpha \not| \beta$),
to bring the entry in the positions $\beta$
to the interval $1,\dots,\alpha -1$.
This completes the description of (RED2).

The steps (RED3) and (RED4) are mirror images of (RED1) and (RED2).
They are performed in a similar way.

\medskip

{\bf Part 2.}

Now we show how we can use the above reduction steps
to get the $\gcd ({\bx})$
of a non-zero element
$$
{\bx}=\mel\alpha\beta A B,
$$

at the position $\alpha$ or $\beta$.

In the loop described below
we will apply our reductions to the appropriate matrix, but we will
concentrate our attention on either the pair in the right column or
the pair in the left column of the matrix $\mel\alpha\beta AB$.

Without loss of generality we can assume that either of $\alpha$ or $B$ is
non-zero (otherwise we apply $\tau$ and the requirement is satisfied).
\medskip

\hspace{10mm}
\begin{minipage}{135mm}
\renewcommand{\baselinestretch}{1}
\renewcommand{\arraystretch}{1}
\normalsize

{\sc Loop.}

{\em Step 1.}
We have an element of the form
$$
\mel{\alpha'} {\beta'} {A'} {B'}
$$
with non-zero pair $\alpha', B'$ in the left column.

By (RED1) we can bring it to the form
$$
\mel {g_1}\beta A 0
$$
with some (new) values of $\beta, A$ and
$0<g_1\le \min\{|\alpha'|, \gcd (B')\}$. It follows from (\ref{end2})
below that $g_1$ is strictly smaller than
the value of $g_2$ from the previous step
({\em ignore this remark if you just entered the loop}).

We need to consider two cases:
\begin{itemize}

\item
$g_1 | \beta, g_1 |A$

We are done: leave the loop.

\item
$g_1$ does not divide either $\beta$ or $A$

It follows that $g_1 \not | \gcd (\beta, A)$.
The assumptions of (RED2) are now satisfied, and we apply it to
$\mel {g_1} \beta A 0$. The result is a new element of the form
$\mel *{\beta''}{A''} {{}*}$ with
\begin{equation} \label{end1}
0<\min\{ \beta'',\gcd (A'')\}<g_1.
\end{equation}

\end{itemize}
Proceed to {\em Step 2}.
\vspace{3mm}

\end{minipage}

\hspace{10mm}
\begin{minipage}{135mm}
\renewcommand{\baselinestretch}{1}
\renewcommand{\arraystretch}{1}
\normalsize

{\em Step 2.}

This step is a mirror image of the previous step, with the roles of elements in the
left column and the right column exchanged.

We have an element of the form
$$
\mel{\alpha''} {\beta''} {A''} {B''}
$$
with elements in the right column satisfying (\ref{end1}).

Using (RED3) we can bring it to the form
$$
\mel\alpha{g_2} 0 B
$$
with some (new) values of $\alpha, B$ and
$0<g_2\le \min\{|\beta''|,\ \gcd(A'')\}$. And (\ref{end1}) implies that~
$g_2$ is strictly smaller than the value of $g_1$ from the previous step.

We consider two cases again:
\begin{itemize}

\item
$g_2 |\alpha, g_2 | \gcd (B)$

We are done: leave the loop.

\item
$g_2$ does not divide either $\alpha$ or $B$

It follows that $g_2 {}\not\!|{\ } \gcd (\alpha, B)$. %
The assumptions of (RED4) are now satisfied, and we apply it to
$\mel \alpha {g_2} 0 B$. The result is a new element of the form
$\mel {\alpha'} **{B'}$ with
\begin{equation} \label{end2}
0<\min\{ \alpha',\gcd (B')\}<g_2.
\end{equation}
\end{itemize}
Proceed to {\em Step 1}.
\end{minipage}

\vspace{4mm}
As we run this loop the values of $g_1$ and $g_2$ will decrease
at each even and odd step, respectively.
On the other hand, they both remain positive
integers, and this guarantees us that we leave the loop at some step.

If we left the loop at the first step, then the first entry of the
resulting element is $\gcd({\bx})$. And if we left it at the second
step, then the $\gcd({\bx})$ is at the second position. In the latter
case we apply $\tau^{-1}$ to bring it to the first entry.

Let us change the notation again (we have done it many times recently),
and denote the obtained element by
$$
\mel\alpha\beta A0
$$
(it follows from the above that the ``$2,1$"-entry is equal to zero).

We can apply a transformation from $\npg\jz$,
if necessary, to convert $A$ to the Smith normal form. Then $A$
turns into a diagonal matrix, and the last requirement of the lemma is
satisfied.

This element has the desired form.
\bigskip

{\bf (b)}
We now consider the case $\jz=\zet\oplus\zet\oplus\zet$
(and hence $\mz=\m(\zet\oplus\zet\oplus\zet)$).
In this case our ``matrices" in $\jz$ are already in the ``diagonal" form.
But the norm-preserving group is too small in this case,
and it is not true that every element of $\jz$ may be converted to a Smith
normal form. However, one can still use the ideas from the algorithm in part (a)
with minor modifications to prove the assertion of the lemma.

We will not pursue the route described in the previous paragraph to prove~(b).
Instead we will use the relation between the pairs
$\Bigl(\grz, \mz\Bigr)$ and
$\Bigl(\SL_2(\zet)\times \SL_2(\zet)\times \SL_2(\zet),\
\zet^2\tensor \zet^2\tensor \zet^2\Bigr)$ described in Example~$\ref{ex:FFF}$.

The reduction procedure described in \cite[Appendix]{bh1}
explains how an arbitrary element of $\zet^2\tensor \zet^2\tensor \zet^2$
may be transformed to a certain $5$-parameter form by a transformation in
the group $\SL_2(\zet)\times \SL_2(\zet)\times \SL_2(\zet)$.
Using the isomorphisms of Example~$\ref{ex:FFF}$,
this procedure may be restated in terms of
transformations in $\grz$ acting on the elements on $\m(\jz)$,
which yields the desired result.
\ep

An immediate corollary of this reduction lemma is Proposition~$\ref{genz}$
that describes the generators of the group $\grz$.

\begin{prop} \label{genz}
Let $\jz$ be one of $\h{\binz}, \h{\quatz}, \h{\octz},
\zet\oplus\zet\oplus\zet$, and let
$\mz=\m(\jz)$.
The group $\grz$ is generated by the elementary transformations
\begin{equation} \label{ele}
\phi(C),\ \psi(D),\ \tf{s},\ \tau
\end{equation}
with $C,D\in \jz,\ s\in \str\jz$.
\end{prop}

{\sc Proof.}

Let $\rho$ be an arbitrary element in the group $\grz$. We need to prove
that $\rho$ can be represented as a product of elementary transformations
$(\ref{ele})$.

Let $f_1$ be
the element $\mel 1000\in \mz$. Let us consider the element
$$
\bx \defeq \rho (f_1).
$$

By Lemma~$\ref{bl-fz}$ we can use a product $\eta$ of elementary
transformations to bring $\bx$ to the (diagonal) reduced form
$$
\bx'=\mel\alpha\beta A 0,\qquad \alpha>0.
$$

We have $\bx' = \eta\, \rho (f_1)$, and it follows from %
Lemma~$\ref{lem:rank}$ that
$$
{\rm rank\ }\bx'={\rm rank\ }f_1=1.
$$

Then it follows from the relations in Remark~$\ref{rem1}$
that $\beta = 0$, $A=0$, and hence
$$
\bx'=\mel{\alpha}000.
$$

We have $\eta, \rho\in \grz$, and it follows from Lemma~$\ref{gcd}$ that
$\gcd(\bx')=\gcd(f_1)$,
and hence we get that
$$
\alpha=1 \qquad {\rm and} \qquad \bx'=f_1,
$$

and hence
$$
\eta \rho (f_1)=f_1.
$$

Repeating the argument of \cite[Lemma 12,13]{b} and using the fact
that our transformations preserve the integral structure,
we prove that $\eta \rho$ has the form
$$
\eta \rho = \tf{s}\phi(C) \quad
{\rm for\ some \ }s\in \str\jz, C\in \jz,
$$
and since $\eta$ was the product of elementary transformations, $\rho$ has
the desired form.
\ep

\subsection{Projective elements in the module $\m(\jz)$}
\label{ss:proj}

In this subsection we discuss the concept of projective elements,
which plays the central role in the classification of orbits in the
integral case.

\begin{defn}\label{def:proja}{\bf (a)}

{\rm
Let $\jz=
\zet\oplus\zet\oplus\zet$, and let $\mz=\m(\jz)$.
An element $x\in \mz$
$$
x=
\mel\alpha\beta
{({a_1}, {a_2}, {a_3})}
{({b_1}, {b_2}, {b_3})},\qquad
\alpha, \beta, a_i, b_i\in \zet
$$
is said to be {\em projective}, if each of the three binary quadratic forms
associated to this element
\begin{eqnarray*}
\biggl(\alpha a_1-b_2 b_3\biggr)x^2 +
\biggl(-a_1b_1+a_2 b_2+ a_3 b_3 -\alpha\beta\biggr)xy +
\biggl(\beta b_1-a_2a_3\biggr) y^2,\\[2mm]
\biggl(\alpha a_2-b_3 b_1\biggr)x^2 +
\biggl(\phantom{-}a_1b_1-a_2 b_2+ a_3 b_3 -\alpha\beta\biggr)xy +
\biggl(\beta b_2-a_3a_1\biggr) y^2,\\[2mm]
\biggl(\alpha a_3-b_1 b_2\biggr)x^2 +
\biggl(\phantom{-}a_1b_1+a_2 b_2- a_3 b_3 -\alpha\beta\biggr)xy +
\biggl(\beta b_3-a_1a_2\biggr) y^2.
\end{eqnarray*}
is primitive, i.e., the $\gcd$ of its coefficients is equal to one.

}
\end{defn}

The concept of projective element was first introduced by M.~Bhargava
\cite{bh-phd, bh1} for elements of
$\zet^2\tensor\zet^2\tensor\zet^2$. The equivalence of the two definitions
follows from the isomorphism
$$
\zet^2\tensor\zet^2\tensor\zet^2 \cong \m(\zet\oplus\zet\oplus\zet)
$$
established in Example~\ref{ex:FFF}.

The concept of projective element was extended to
$\wedge^3\zet^6$ (and several more examples) in~\cite{bh1}.
More precisely, M.~Bhargava showed that there is a natural injection
\begin{equation}\label{inc1}
\zet^2\tensor\zet^2\tensor\zet^2 \hookrightarrow \wedge^3\zet^6
\end{equation}
and in addition every element of $\wedge^3\zet^6$ is
$\SL_6(\zet)$-equivalent to an element in the image of
$\zet^2\tensor\zet^2\tensor\zet^2$ under the embedding (\ref{inc1}).
An element of $\wedge^3\zet^6$ was defined to be {\em projective},
if its $\SL_6$-orbit contains an image under (\ref{inc1}) of a
projective element in $\zet^2\tensor\zet^2\tensor\zet^2$.

In our case we have a chain of natural inclusions of cubic vector spaces
($\zet$-modules):
$$
\zet\oplus\zet\oplus\zet \subset
\h\binz \subset \h\quatz \subset \h\octz
$$
(the first inclusion being the diagonal embedding), which induces the
natural chain
\begin{equation} \label{chain3}
\m(\zet\oplus\zet\oplus\zet) \subset
\m(\h\binz) \subset \m(\h\quatz) \subset \m(\h\octz).
\end{equation}

The first two modules in the last chain are isomorphic to
$\zet^2\tensor\zet^2\tensor\zet^2$ and $\wedge^3\zet^6$, respectively
(see Examples~$\ref{ex:FFF}, \ref{ex:sl6int}$). We proved in Lemma~\ref{bl-fz}
that an arbitrary element of
each of the modules $\mz$ in (\ref{chain3}) is $\grz$-equivalent to
a {\em diagonal reduced} element,
which can be thought of as an element of the smallest submodule
$\m(\zet\oplus\zet\oplus\zet)$ in the chain~$(\ref{chain3})$.
This lemma allows us to extend
the concept of a projective element to $\m(\h\binz)\cong \wedge^3\zet^6$,
using an argument different from Bhargava's, and further extend this
concept to $\m(\h\quatz)$ and  $\m(\h\octz)$.

\hspace{-\parindent}{\bf Definition~\ref{def:proja} (b)}
{\rm
Let $\jz$ be one of $\h{\binz}, \h{\quatz}, \h{\octz}$, and let
$\mz=\m(\jz)$. An element $x\in\mz$ is said to be {\em projective}, if its
$\grz$-orbit contains a diagonal reduced element (cf. Lemma~$\ref{bl-fz}$)
\begin{equation}\label{proj-el}
\mel\alpha\beta{\diag{a_1} {a_2} {a_3}}{0},\qquad
\alpha, \beta, a_i\in \zet,\quad
\alpha>0,\quad \alpha |\beta, \ \alpha |a_i,
\end{equation}
which is projective.
}

We note that for an element of the form (\ref{proj-el}) the projectivity
conditions of Definition~\ref{def:proja}(a) become
\begin{eqnarray} \label{def:projexpr}
& \gcd(\alpha a_1, \alpha\beta, a_2a_3)=1 & \nonumber\\
& \gcd(\alpha a_2, \alpha\beta, a_1a_3)=1 &\\
& \gcd(\alpha a_3, \alpha\beta, a_1a_2)=1 & \nonumber
\end{eqnarray}
Taking into account the divisibility conditions of reduced elements
(Definition~\ref{def:reduced}), it follows that if element (\ref{proj-el})
is projective, then $\alpha=1$, i.e., $x$ is primitive. Nevertheless we
kept $\alpha$ in the expressions (\ref{def:projexpr}) to emphasize the fact
that they are homogeneous expressions of degree $2$.
\medskip

The definition of a projective element we just gave is not quite
satisfactory, since
it is not convenient to work with orbit representatives of the group
$\grz$, which may be quite large.
Next we would like to argue that the concept of a projective
element is related to certain equations of degree $3$, which will allow us
to get a simple projectivity test in {\em almost} all cases.

\begin{prop} \label{prop:proj}
{\rm
Let $\jz$ be as in Definition~$\ref{def:proja}$(a,b) and $\mz=\m(\jz)$.
Let $x$
be in $\mz$ and let $T(x,x,x)$ be defined as in~$(\ref{deft})$.
\begin{enumerate}[(i)]
\item
If $\gcd\, T(x,x,x)=1$ then $x$ is projective;
\item
If $\gcd\, T(x,x,x)\ge 3$ or $T(x,x,x)=0$ then $x$ is not projective.
\item
When $q'(x)$ is odd, $\gcd T(x,x,x)=1$ iff $x$ is projective.
\end{enumerate}
}
\end{prop}

{\sc Proof.}

First we notice that the quantity $\gcd T(x,x,x)$ is invariant with
respect to the action of the group $\grz$.
This assertion follows from the
$\gr$-invariance of the quartic form $q(x)$ and the symplectic form
$\{\cdot,\cdot\}$,
the definition of $T(x,x,x)$ $(\ref{deft})$, and Lemma~$\ref{gcd}$
(see Lemma~$\ref{di}$ below for a detailed argument).
Hence, for a diagonal reduced element
\begin{equation}\label{x1a}
x_1=\mel\alpha\beta{\diag{a_1} {a_2} {a_3}}{0},\qquad
\alpha, \beta, a_i\in \zet,\quad
\alpha>0, \alpha|\beta, \alpha| a_i.
\end{equation}
contained in the $\grz$-orbit of $x$ by Lemma~$\ref{bl-fz}$ we have
$$
\gcd T(x,x,x)=\gcd T(x_1,x_1,x_1).
$$

It is more convenient to work with $x_1$,
since the components
of $T(x_1,x_1, x_1)$ have the following simple form (see~$(\ref{txxx})$):
$$
T(x_1,x_1,x_1)=
\mel{\alpha^2\beta}{\ \alpha\beta^2+2\mnorm(A)}{\ \alpha\beta A}{\ 2\alpha
A^\#},
$$

where $A=\diag{a_1} {a_2} {a_3}$,
$\mnorm (A)=a_1a_2a_3,\ A^\#=\diag{a_2a_3}{a_3a_1}{a_1a_2}$.

{\bf (i)}
Let $\gcd T(x_1,x_1,x_1)=1$, and suppose $x_1$ is not projective. It means
that the $\gcd$ of at least one of the expressions in~$(\ref{def:projexpr})$
is greater than $1$. Without loss of generality we assume that there
exists a prime $p>1$ such that
$$
p\ \Bigl|\  \alpha a_3,\  \alpha\beta,\ a_1a_2.
$$
This implies that $p$ divides each component of $T(x_1, x_1, x_1)$,
contrary to our assumption.

{\bf (ii)}
We are given $x$ such that $\gcd T(x, x, x)\ge 3$, and need to show that
$x$ is not projective. Let $x_1$ be {\em any} diagonal reduced element of
the form~$(\ref{x1a})$ in the $\grz$-orbit of $x$. Hence
$$
\gcd T(x_1, x_1, x_1)=
\gcd
\mel{\alpha^2\beta}{\alpha\beta^2+2\mnorm(A)}{\alpha\beta A}{2\alpha A^\#}
\ge 3.
$$

First we notice that if $\alpha >1$, then $x_1$ is not primitive.
Non-primitivity implies non-projectivity, and hence there is nothing to
prove in this case.
Thus from now on we may assume
$\alpha=1$, and hence
$$
T(x_1, x_1, x_1)=
\mel{\beta}{\beta^2+2\mnorm(A)}{\beta A}{2 A^\#}.
$$

Suppose that there exists a prime $p\ge 3$ that divides $\gcd T(x_1,
x_1, x_1)$. Then we have $p|\beta$. In addition, the condition $p|A^\#$
implies that at least two of the numbers $a_i$ are divisible by $p$. Then
$p$ divides at least one of the expressions~$(\ref{def:projexpr})$, and
hence $x_1$ is not projective.

If the prime $p$ such as in the previous paragraph does not exist, it
follows from the conditions that $\gcd T(x_1, x_1, x_1)$ is a power of $2$
(and is greater than $2$). Then we have $4|\beta$, $2|A^\#$, and repeating
the argument of the previous paragraph, we get that $2$ divides at least
one of the expressions~$(\ref{def:projexpr})$. Hence $x_1$ is not
projective in this case either.

We have proved that if $\gcd T(x,x, x)\ge 3$, then every diagonal reduced
element in $\grz$-orbit of $x$ is not projective, and hence $x$ is not
projective.

This argument also applies in the case $T(x,x,x)=0$, if one replaces the
condition ``divisible by $p$" with the condition ``equal to $0$".

{\bf (iii)}
The implication ``$\Rightarrow$" was proved in (i). To prove the converse
implication, we again look at the diagonal reduced projective
$x_1$ of the form $(\ref{x1a})$ in the $\grz$-orbit of $x$.

It follows from (ii) that $\gcd T(x_1,x_1,x_1)$ is equal to $1$ or $2$,
and we will show that the latter option is not possible.

Since $x_1$ is projective, it is primitive, and hence $\alpha=1$.
We then have $q'(x_1)=\beta^2+4\mnorm (A)$, and since $q'(x)=q'(x_1)$ is
odd, we have that $\beta$ is odd. This conclusion implies that $2$ cannot
divide $\gcd T(x_1,x_1,x_1)$, and hence $\gcd T(x_1,x_1,x_1)=\gcd T(x,x,x)=1$.
\ep

\begin{rem}
{\rm
The above proposition  provides a projectivity test for all cases
except $\gcd\,$ %
$T(x,x,x)=2$. Elements
\begin{equation}
x_1=\mel{1}{2}{\diag{1}12}0 \qquad
x_2=\mel{1}{2}{\diag{1}22}0
\end{equation}
both have $\gcd \Bigl(T(x_i,x_i,x_i)\Bigr)=2$.
However $x_1$ is projective and $x_2$ is not projective in the sense of
Definition~$\ref{def:proja}(a)$.
}
\end{rem}

\begin{cor}\label{cor:anyrep}
{\rm
Proposition~$\ref{prop:proj}$ implies that if $\gcd T(x,x,x)$ is
different from $2$ then the (non)proje\-ctivity of $x$ is determined by
{\em any} representative of its $\grz$-orbit
in $\m(\zet\oplus\zet\oplus\zet)$.
}
\end{cor}

It appears that the assertion of Corollary~$\ref{cor:anyrep}$
remains true for $x$ with $\gcd\, T(x,x,x)=2$,
however the quantity $T(x,x,x)$ does not seem to be precise enough to
treat this case. Apparently, the difficulty of proving this fact
is related to the ``bad" reduction
of the quartic form in the case when char\,$F=2$ (cf. Remark~$\ref{badred}$).
\medskip

Our next remark is concerned with one more characterization of projective
elements.
There is a natural choice of a basis and coordinates in the module
$\mz=\m(\jz)$ with $\jz=\zet\oplus\zet\oplus\zet,\h\binz,\h\quatz,\h\octz$.
In each of these cases the quartic form $q'(x)$ may be viewed as a
homogeneous polynomial in $n$ variables, $n=8, 20, 32, 56$, respectively.
And the components of $T(x,x,x)$ are nothing else but the $n$ partial
derivatives of the polynomial $q'(x)$ (divided by $2$, since all
these partial derivatives have even coefficients).
In other words, each component of the formal gradient vector
$\frac{1}{2}\nabla q'(x)$ is a polynomial with integer coefficients,
and these expressions are the same as the components of $T(x,x,x)$.
These observations and Proposition~$\ref{prop:proj}$ yield the following

\begin{cor}\label{cor:nabla}
{\rm
Let $x$ be in $\mz=\m(\jz)$ with $\jz$ as in Definition~$\ref{def:proja}$(a,b).
Let $q'$ be the quartic invariant of the module $\mz$, see~$(\ref{qmod})$.
Then $\gcd T(x,x,x)=\gcd \frac{1}{2}\nabla q'(x)$ and hence
\begin{itemize}
\item
If $\gcd\, \frac{1}{2}\nabla q'(x)=1$ then $x$ is projective;
\item
If $\Bigl(\gcd\, \frac{1}{2}\nabla q'(x)\ge 3\Bigr)$
or $\Bigl(\nabla q'(x)=0\Bigr)$ then $x$ is not projective.
\end{itemize}
}
\end{cor}

\begin{rem}\label{rem:nabla}
{\rm
It follows
from general considerations that if $p(x)$ is a (quartic) form on a
$\zet$-module $V_\zet$ invariant with respect to a group $G_\zet$, then
$\gcd \nabla p(x)$ is invariant with respect to $G_\zet$. This observation
suggests that Corollary~$\ref{cor:nabla}$ may be used to describe
(non)projective elements in spaces ${\rm Sym}^3\zet^2,
\zet^2\tensor{\rm Sym}^2\zet^2, \zet^2\tensor\zet^4$
underlying higher composition laws
related to quadratic rings (see Table~1 and \cite{bh1} for details).
}
\end{rem}

\subsection{Further reduction and the classification of the projective orbits}
\label{ss:fred}

\begin{lemma}[\sc Reduction Lemma II]  \label{red2}
\

\hspace{-\parindent}%
Let $\jz$ be one of $\h{\binz}, \h{\quatz}, \h{\octz}$
and $\mz=\m(\jz)$.
Let ${\bx}\in \mz$ be a projective element.

Then $x$ is equivalent to an element
\begin{equation} \label{proj1}
\mel 1 \varepsilon{\ \diag{1}1k} 0,\qquad
\varepsilon\in\{0,1\},\quad
k\in \zet
\end{equation}

under a series of elementary transformations
\begin{equation} \label{elem-tr1}
\phi (C), \psi(D), \tf{s}, \tau \quad
with \ C, D \in \jz,\ s\in \npg{\jz}.
\end{equation}

The values of $\varepsilon$ and $k$ in~$(\ref{proj1})$ are uniquely
determined by $q'(x)$.
\end{lemma}

{\sc Proof.}

By the definition of a projective element, $x$ is equivalent to a
diagonal reduced projective element
$$
x_1=\mel{\alpha_1} {\beta_1}{\ A_1} 0, \qquad
{\rm where}\ \alpha_1=\gcd(x_1)=\gcd(x).
$$
under the action of transformations (\ref{elem-tr1}). Since any projective
element is primitive, we have $\alpha_1=1$. In addition, acting by $T(s)$
if necessary we can assume that $A_1$ is in the Smith normal form
(Theorem~$\ref{np-orbits}$).
Summarizing these remarks we conclude that $x_1$ is of the form
\begin{equation}\label{x1}
x_1=\mel 1{\beta_1}{\ \diag{a_1}{a_2}{a_3}}0
\end{equation}
with $\diag{a_1}{a_2}{a_3}$ in the Smith normal form.

The definition of the projective element
(Definition~$\ref{def:proja}$(a), see also~$(\ref{def:projexpr})$) implies that
\begin{equation} \label{proj3}
\gcd\{ \beta_1,\ a_1 a_2, \ a_1 a_3, \ a_2 a_3\}=1.
\end{equation}
Since $A_1$ is in the Smith normal form with $a_1 | a_2,\ a_2|a_3$,
the relation (\ref{proj3}) is equivalent to
\begin{equation} \label{proj4}
\gcd\{ \beta_1,\ a_1 a_2\}=1.
\end{equation}

{\bf Step 1.}
We show that $a_1$ can be taken to be equal to $1$.

Assume this is not the case. Then $a_1>1$
and $a_1|a_2, a_1|a_3$ (case $a_1=a_2=a_3=0$ needs slightly
different treatment; we skip details).
Relation (\ref{proj4}) implies $\gcd(\beta_1, a_1)=1$.

Then we apply Lemma~\ref{lcomp} (iii) to $x_1$ (with $c=1$)
and transform it to
\begin{equation}\label{proj5}
\mel 1
{\beta_1-2a_1a_2}{\quad \diag{a_1}{\ a_2}{\ a_3+\beta_1-{a_1a_2}}}0,
\end{equation}
It follows from the above that
$\gcd\{{a_1},{a_2},\ a_3+\beta_1-{a_1a_2}\}=
\gcd\{{a_1},{a_2}, \beta_1\}=1$.
The $\gcd$ condition implies that the Smith canonical form
of $$\diag{a_1}{a_2}{\ a_3+\beta_1-{a_1a_2}}$$ looks like
$\diag{1}**$.
Then we can apply an appropriate $T(s)$
to (\ref{proj5}),
and it will yield the desired result, completing the first step.

\medskip
Step 1 implies that we can assume that $a_1$ in the element~$(\ref{x1})$
is equal to $1$ and we proceed to

{\bf Step 2.}
We show that $a_2$ in~$(\ref{x1})$ can be taken to be equal to $1$.

Assume this is not the case. Then $a_2>1$ and we still have $a_2|a_3$
(case $a_2=a_3=0$ treated similarly).
Relation (\ref{proj3}) implies $\gcd(\beta_1, a_2)=1$.

We again apply Lemma~\ref{lcomp} (iii) to $x_1$ (with $c=1$)
and transform it to
\begin{equation}\label{proj6}
\mel 1
{\beta_1-2a_2}{\quad \diag{1}{\ a_2}{\ a_3+\beta_1-{a_2}}} 0.
\end{equation}

We have $\gcd\{ {a_2}, a_3+\beta_1-{a_2} \}=
\gcd\{ {a_2},\beta_1\}=1$.
Similarly to Step~1, we can apply an appropriate $T(s)$
to (\ref{proj6}), an get an element of the form
$$
\mel 1{\beta_1-2a_2}{\quad \diag{1}{1}{*}}0.
$$
The second step is complete.

It follows from the above steps that element (\ref{x1}) may be taken to be
of the form
\begin{equation}\label{proj7}
\mel 1 {\beta_1}{\quad \diag{1}{1}{a_3}}0.
\end{equation}

{\bf Step 3.}
We show that $\beta_1$ in the element (\ref{proj7}) may be taken to be $0$ or
$1$.

To do it we again apply Lemma~\ref{lcomp} (iii):
$$
\mel 1{\beta_1-2c}{\quad \diag{1}{\ 1}{\ a_3+\beta_1 c-c^2}}0.
$$

Obviously, there is a $c$, which makes the second component $0$ or
$1$.
\bigskip

The three steps above show that an arbitrary projective element can be
brought to the form~(\ref{proj1}). It remains to prove the uniqueness
assertion.

When $x'$ is of the form (\ref{proj1}), we obtain from
$(\ref{qmod})$ that
$$
q'(x') = 4 k+\varepsilon ^2 .
$$

This equation determines the parity of $\varepsilon$ uniquely,
and hence $k$ is unique as well.
\ep

Now we have all the necessary tools to prove the main result of the paper:
the theorem concerning orbits in the spaces associated to cubic Jordan
algebras of Hermitian matrices over split composition algebras.

\begin{theorem}\label{thm:main}

Let $(G_\zet, \mz)$ be one the following pairs
$$
\Bigl(\SL_6(\zet), \wedge^3(\zet^6)\Bigr),\quad
\Bigl(D_6(\zet), \mbox{\rm half-spin}_\zet\Bigr),\quad
\Bigl(E_7(\zet), V(\omega_7)_\zet\Bigr).
$$
Then
\begin{itemize}
\item
The $G_\zet$-invariant quartic form (the norm) on the module $\mz$
has values congruent to $0$ or $1\, (\!\mod 4)$.

\item
Let $n$ be an integer $\equiv 0$ or $1\, (\!\mod 4)$.
The group $G_\zet$ acts
transitively on the set of {\em projective} elements of norm $n$.

\item
If $n$ is a fundamental discriminant\footnote{
An integer $n$ is called a fundamental discriminant if $n$ is squarefree
and $\equiv 1 (\mod 4)$
or $n=4k$, where $k$ is a squarefree integer that is $\equiv 2$ or $3 (\mod 4)$.
The result stated in the theorem also applies in the case $n=1$.
},
then every element of norm $n$ is projective, and hence in this case
$G_\zet$ acts transitively on the set of elements of norm $n$.
\end{itemize}
\end{theorem}

{\sc Proof.}
It was shown in Proposition~$\ref{prop:gener}$ that for $\J=\h\bin,
\h\quat, \h\oct$ the Freudenthal construction yields an
absolutely almost simple connected algebraic
group $G=\gr$ of type $A_5, D_6, E_7$, respectively.
Each of these groups acts on the vector space $\m=\mj$, producing
an irreducible representation
whose type is given in the statement of the theorem.

The integral structure in the group $G_\zet=\grz$
and the module $\mz$ is induced by the integral structure in $\jz$.
The case $\jz=\h\binz$ requires special treatment.
It was shown in Example~$\ref{ex:sl6int}$ that $\SL_6(\zet)$ is isomorphic
to a subgroup (of index two) in $\g{\m(\h{\bin_\zet})}$. We will do the proof
of the theorem for the whole group $\grz$, and address the issue of
$\SL_6(\zet)$ at the end of the proof.

The module $\mz=\m(\jz)$ comes equipped with the quartic form
$$
q'(x)= \Bigl((A,B)-\alpha\beta\Bigr)^2
-4 (A^\#, B^\#) +4\alpha \mnorm(A) +4\beta \mnorm(B),\qquad
x=\mel{\alpha}\beta A B \in \mz,
$$
invariant with respect to $G_\zet$. It was noted in Subsection~$\ref{ss:33int}$
that the operations $ \mnorm$ and $(\cdot,\cdot)$ in $\jz$ have integer
values, which implies that $q'(x)$ is always congruent
to $0$ or $1$ modulo $4$.

Next, let $x$ be a projective element in $\mz$. It follows from
Lemma~$\ref{red2}$, that $x$ can be transformed to an element of the
form
\begin{equation}
\mel{1}\varepsilon {\diag{1}1k}0,\qquad
\varepsilon\in\{0,1\},\quad
k\in \zet
\end{equation}
with values of $\varepsilon$ and $k$ uniquely determined by the value of
$q'(x)$. This immediately implies that the group $\grz$ acts transitively
on the set of projective elements of norm $q'(x)$.

Finally, let an integer $n$ be a fundamental discriminant, and let
$x\in\mz$ be such that $q'(x)=n$.
By Lemma~$\ref{bl-fz}$, $x$ is equivalent to a diagonal reduced element of the
form
\begin{equation}\label{red-thm}
x_1=\mel \alpha\beta {\diag{a_1}{a_2}{a_3}} 0,\qquad
\alpha>0,\alpha|\beta,\alpha|a_i.
\end{equation}
We need to show that $x$ is projective, and by
Definition~$\ref{def:proja}$(b), it is sufficient to show that $x_1$ is
projective.

We have $q'(x)=q'(x_1)$, and hence
\begin{equation}\label{n1}
n=\alpha^2\beta^2+4a_1a_2a_3.
\end{equation}

First we show that $\alpha$ must be equal to $1$.

Note that $(\ref{red-thm})$ and $(\ref{n1})$ imply $\alpha^3|n$.

If $\alpha >2$, this remark implies that $n$ is not square-free. If
$\alpha=2$, we get $16|n$. And in both cases we get that $n$ is not a
fundamental discriminant. Hence $\alpha=1$.

It means we can rewrite $x_1$ in the form
$x_1=\mel 1 \beta {\diag{a_1}{a_2}{a_3}} 0$, and
$$
n=\beta^2+4a_1a_2a_3.
$$

To complete the proof we will show that if $x_1$ is not projective, then
$n$ is not a fundamental discriminant.

So suppose that $x_1$ is not projective. Then the $\gcd$ in one of the
expressions~$(\ref{def:projexpr})$ is greater than $1$. Without loss of
generality (and using $\alpha=1$) we will assume that
$$
\gcd (a_3, \beta, a_1 a_2)>1.
$$
Let $p$ be a prime dividing $\gcd (a_3, \beta, a_1 a_2)$. It follows that
$$
p|\beta \quad \mbox{ and \ $p$ divides at least two of the $a_i$'s}.
$$
This remark implies $p^2|n$. If $p>2$ it already implies that $n$ is not a
fundamental discriminant. And if $p=2$, it implies that
$$
\beta=2\beta_1, \qquad
a_1a_2a_2=4c
$$
for some integers $\beta_1$ and $c$. Hence $n$ may be rewritten in the
form
$$
n=4\Bigl(\beta_1^2+4c\Bigr).
$$
Depending on the parity of $\beta_1$, the quantity $\beta_1^2+4c$ is
congruent to $0$ or $1$ modulo $4$. In both cases it implies that $n$ is
not a fundamental discriminant.

We thus proved that if $n$ is a fundamental discriminant, then every
element of norm $n$ is projective, and hence $G_\zet$ acts transitively on
the set of such elements.
\medskip

Our last remark is concerned with the case $\jz=\h\binz$. It was noted
in Example~$\ref{ex:sl6int}$ that $\SL_6(\zet)\cong \g{\jz}\cc$, which is
a subgroup of index two in $\grz$. We need to show that orbits under
$\g{\jz}\cc$ are the same as the orbits under the action of the whole
group $G_\zet=\grz$.

For this we take an arbitrary element $x\in\mz$ and let $x_1$ be a
diagonal reduced element
$$
x_1=\mel \alpha\beta {\diag{a_1}{a_2}{a_3}} 0,\qquad
\alpha>0,\alpha|\beta,\alpha|a_i
$$
lying in the $\grz$-orbit of $x$, i.e., $x_1=g(x)$ for some $g\in\grz$.
It was noted in Example~$\ref{ex:sl6int}$ that every $g\in\grz$ has the
form
$$
g=g_0\quad
\mbox{or}\quad
g=T(t) g_0
$$
with $g_0\in\grz\cc$. There is nothing to prove in first case, and in the
second case we have
$$
x_1=T(t)(x_1)=T(t)^2\,g_0 (x)=g_0(x)
$$
using that $T(t)$ acts as the transpose operation
on the matrices at the two off-diagonal
entries of $x_1$, and the fact that $T(t)^2=\id_\mz$. Hence $\grz$-orbits
are the same as $\grz\cc$-orbits when $\jz=\h\binz$, and in this case
the theorem remains true for the subgroup of index two isomorphic to
$\SL_6(\zet)$.
\ep

\subsection{Invariant factors
in the module $\mz$ and the degenerate orbits}
\label{ss:if}

In this subsection we introduce the ``invariants" $d_i, i=1,2,3,4$,
reminiscent of invariant factors
for regular matrices over integers, and use them to
describe the orbits of degenerate elements, i.e., those whose norm is
equal to zero.
We define these invariants using the concept of the rank polynomials introduced
earlier in $(\ref{p4})-(\ref{p2})$.

\medskip

\begin{defn} \label{def:di}
{\rm Let $\jz$ be as in~$(\ref{zlist})$, and let $\mz=\m(\jz)$.
For an element $\bx\in\mz$ of the form
$$
\bx = \mel\alpha\beta A B
$$
the functions $d_i:\ \mz\ \to \ \zet$ are defined by the following
expressions
\begin{itemize}
\item
$d_1 (\bx) =\gcd(\bx)$;
\item
$d_2 (\bx) = \gcd \Bigl(3\,T(x,x,y)+\{x,y\}\,x\Bigr) \quad
\mbox{for all }y\in \mz$;
\item
$d_3 (\bx) = \gcd ( T(\bx,\bx,\bx) )$;
\item
$d_4 (\bx) = q'(\bx)$.
\end{itemize}
}
\end{defn}

\begin{lemma} \label{di}
Let $\mz$ be as in Definition~$\ref{def:di}$.
\begin{enumerate}[{\rm (i)}]
\item
The functions $d_i$ are invariant under the action of the
group $\grz$.
\item
\begin{equation}\label{d2}
d_2(x)= \gcd\Bigl\{
3\alpha\beta-(A,B),\quad
2(\alpha A-B^\# ), \quad 2(\beta B - A^\# ),\quad
2Q(x)\Bigr\}.
\end{equation}
\end{enumerate}

\end{lemma}

{\sc Proof.}

{\bf (i)}
The statement for $d_1$ follows from Lemma~$\ref{gcd}$.

The function $d_4$ is the norm $q'$, and the statement follows from the
definition of the groups $\gr$ and $\grz$.

The statement for $d_3$ follows from the following computation
$$
d_3(\eta(x))= d_1\biggl(\  T\Bigl(\eta(x), \eta(x),\eta(x)\Bigr)\
\biggr)=
$$
$$
=d_1\biggl(\ \eta \Bigl( T(x,x,x) \Bigr)\ \biggr)= d_1\biggl( T(x,x,x)
\biggr)= d_3(x) \qquad \mbox{for any }\eta\in \grz.
$$

In this computation we used the definition of $d_3$, the relation
$(\ref{autot})$, and the invariance statement for $d_1$.

\medskip
The statement for $d_2$ follows from the invariance of the skew-symmetric
form $\{\cdot, \cdot\}$ and the argument analogous to that in the previous
paragraph.

{\bf (ii)}
The computation ~$(\ref{r1eq})$ of Lemma~$\ref{lem:rank1}$ implies that
$$
\mbox{an integer $d$ divides}\quad
3\,T(x,x,y)+\{x,y\}\,x \quad
\mbox{for any }y\in \mz
$$
if and only if $d$ divides each of the following expressions
$$
3\alpha\beta-(A,B),\quad
2(\alpha A-B^\# ), \quad 2(\beta B - A^\# ),\quad
2Q(x),\quad 2Q(x'),
$$
where $x=\mel \alpha\beta A B$, $x'=\mel \beta \alpha B A$, and
$Q(x)\in \End_\zet(\jz)$ was defined in~$(\ref{q})$.

As in the proof of Lemma~$\ref{lem:rank1}$ we have
$$
Q(x')(C)=Q(x)(C)-2Q(x)(\one)\jprod C\quad
\mbox{for any }C\in \jz.
$$
This implies that $d$ divides $2Q(x')$ whenever $d$ divides $2Q(x)$,
and the statement (ii) of the lemma follows.
\ep

\begin{rem}\label{rem:rpz}
{\rm
It is possible to define alternative quadratic invariant $\dvap$
via
\begin{equation}\label{alt-d2}
\dvap(x)=\gcd\Bigl(
\alpha A-B^\# , \quad \beta B - A^\# ,\quad Q(\bx)
\Bigr).
\end{equation}

However the short proof of $\grz$-invariance of $d_2$ in the previous
lemma does not work for $\dvap$. One can still prove that $\dvap$ is
invariant under each of the transformations $\phi(C), \psi(D), T(s)$, and
then using Proposition~$\ref{genz}$ conclude that $\dvap$ is invariant
under the whole $\grz$. This route is quite technical, it was implemented
in Lemma~$2.3.5$ of \cite{kphd}.

Our primary use for the invariants $d_i$ is to distinguish elements of $\mz$
lying in different orbits.
We should note that $\dvap$ is ``finer" than $d_2$ in this sense.
For example, using $\dvap$ one can conclude that elements
$$
\mel 1 2 {\diag{1}00} 0
\quad {\rm and}\quad \mel 1 2 {\diag{2}00} 0
$$
lie in distinct $\grz$-orbits, though all $d_i$'s are equal for these two
elements.

In this paper will make use of more ``coarse" relations~$(\ref{d2})$,
which are sufficient for our purposes. \ee
}
\end{rem}

\begin{rem} \label{rem2}
{\rm
By analogy with Remark~$\ref{rem1}$ we notice that when $x$ is a reduced
element
$$
x=\mel \alpha\beta A 0
$$
then
$%
d_2(x)=\gcd\Bigl(\alpha\beta,\ 2\alpha A,\ 2A^\# \Bigr)
$. \ee%
}
\end{rem}

The following theorem in an extended version of Theorem~$\ref{thm:main}$.
It provides the description of degenerate $\grz$-orbits (corresponding the
zero value of the quartic form) in the module $\mz$.

\begin{theorem} \label{thm2}
Let $\jz$ be one of $\h{\binz}, \h{\quatz}, \h{\octz}$
and $\mz=\m(\jz)$.

\begin{itemize}

\item
Every element $\bx$ of rank $1$ in the module $\mz$ can be brought to the
form
\begin{equation} \label{thr1}
\mel \alpha 000 \quad {\rm where}\  \alpha = d_1(\bx)
\end{equation}
by an element in the group $\grz$.

The set
\begin{equation} \label{thset1}
\left.\left\{\mel k 0 0 0\ \right|\ k\in\zet, k\ge 1\right\}
\end{equation}
is a complete set of {\em distinct\/} orbit representatives of elements of
rank $1$ in $\mz$.

\item
Every element $\bx$ of rank $2$ in the module $\mz$ can be brought to the
form
\begin{equation} \label{thr2}
\mel \alpha 0 {\ \diag{a}00} 0
\quad {\rm where}\  \alpha|a,\ \alpha =
d_1(\bx), \ a= d_2(\bx)/\alpha
\end{equation}
by an element in the group $\grz$.

The set
\begin{equation} \label{thset2}
\left\{
\left.\mel k 0 {\diag{m}00} 0 \
 \right.|\ k,m\in\zet, k,m>0, k|m\right\}
\end{equation}
is a complete set of {\em distinct\/} orbit representatives of elements of
rank $2$ in $\mz$.

\item
In the case $\rank x=3$ or $4$, the group $\grz$ acts transitively on the
set of {\em projective} elements of norm $n$.
Every such element may brought to the form
\begin{equation}\label{hirank}
\mel 1 \varepsilon{\ \diag{1}1k} 0,\qquad
\varepsilon\in\{0,1\},\quad
k\in \zet,
\end{equation}
where $k=\frac{q'(x)-\varepsilon^2}{4}$ and $\varepsilon\equiv
q'(x)\,(\mod 2)$ are uniquely determined by $x$.
\end{itemize}
\end{theorem}

{\sc Proof.}

We will often use in the proof the facts that transformations in the group
$\grz$ preserve the rank and the invariants $d_i$ (Lemma~$\ref{lem:rank}$,
$\ref{di}$) of elements of $\mz$.

We note that if $\rank x \le 2$, then $T(x,x,x)=0$
(Definition~$\ref{def:rank}$).
Such an $x$ is not projective by Proposition~$\ref{prop:proj}$(ii),
and hence Theorem~$\ref{thm:main}$ gives no information about orbits
of such elements.

\medskip
{\em The case of\/} rank $1$.

Let $\bx$ be an arbitrary element of rank $1$. First we show that $\bx$
can be brought to the form $(\ref{thr1})$. Then we show that elements of
the form $(\ref{thset1})$ with distinct $k$'s lie in distinct orbits of
the group $\grz$.

By Lemma~$\ref{bl-fz}$ there exists $\sigma\in \grz$ such that
$$
\sigma (\bx) = \mel \alpha \beta A 0 ,
\quad \qquad \alpha>0,\
\alpha|\beta,\ \alpha |A.
$$

Since rank $\bx =1$, it follows from $(\ref{r1s})$ that $\beta =0$ and
$A=0$. So
$$
\sigma (\bx) = \mel \alpha 000
$$

is in the desired form.

\medskip
Now suppose we have two elements $\bx_1=\mel {k_1}000$ and
$\bx_2=\mel{k_2}000$
of the form $(\ref{thset1})$ lying in the same $\grz$-orbit. Then
$d_1(\bx _1)=k_1$ and $d_1(\bx _2)=k_2$, and since $d_1$ is constant on
orbits, we have $k_1=k_2$, and hence $x_1=x_2$.
The proof in the case of rank $1$ is complete.

\medskip
{\em The case of\/} rank $2$.

Let $\bx$ be an arbitrary element of rank $2$. First we show that $\bx$
can be brought to the form~$(\ref{thr2})$. Then we show that elements of
the form~$(\ref{thset2})$ with distinct $k$'s and $m$'s lie in distinct
orbits of the group $\grz$.

By Lemma~$\ref{bl-fz}$ there exists $\sigma\in \grz$ such that
$$
\sigma (\bx) = \mel\alpha\beta A 0, \quad \qquad \alpha>0,\
\alpha|\beta,\ \alpha |A.
$$

Since rank $\bx =2$, it follows from $(\ref{r2s})$ that $\beta =0$ and
$A^\# =0$. The last relation implies that rank $A\le 1$. On the other
hand, $A\ne 0$ as $A=0$ would imply that rank $x \le 1$ (see
$(\ref{r1s})$). Hence rank $A=1$.

By Theorem~$\ref{np-orbits}$,
there exists $\xi\in \npg\jz$ which brings $A$ to the
Smith normal form. Since $\rank A =1$,
this normal form has only one nonzero (in fact, positive)
element on the diagonal. We thus get

$$
\tf{\xi}\, \sigma \ (\bx)= \mel \alpha 0 {\diag{a_1}00} 0.
$$

By construction $a_1 = \gcd(A)$, and so $\alpha | a_1$.

Hence we have brought $\bx$ to the desired form~$(\ref{thr2})$.
\medskip

Now let us take two elements $x_1, x_2$ of the form $(\ref{thset2})$
$$
x_i= \mel {k_i} 0 {\diag{m_i}00} 0\ \qquad k_i,m_i\in\zet, k_i, m_i > 0,
k_i|m_i, i=1,2
$$

lying in the same $\grz$-orbit. We want to show that $x_1=x_2$.

We have $d_1(x_i)=k_i$, and using the fact that $d_1$ is constant on
orbits, we get
$$
k_1=d_1(x_1)=d_1(x_2)=k_2.
$$

Next we notice that Remark~$\ref{rem2}$ implies
$$
d_2(x_i)\ =\ 2\, k_i\, m_i
$$

and hence
$$
2\,k_1\, m_1=d_2(x_1)=d_2(x_2)=2\, k_2\, m_2,
$$

which implies that in this case $m_1=m_2$, and hence $x_1=x_2$.
This completes the proof in the case of rank $2$.

\medskip
{\em The case of\/} rank $>2$.
A complete reduction procedure is not known in the case of elements of
rank $3$ and $4$, and the structure of $\grz$-orbits may be quite complicated.
For example, when $\jz=\h\binz$, the structure of orbits
in $\mz$ is essentially equivalent to the structure of $\SL_6(\zet)$-orbits
in $\wedge^3(\zet^6)$ (see Example~$\ref{ex:sl6int}$).
It is as complicated as the structure of
(balanced) triples of ideal classes in quadratic rings
\cite[Theorem~$18$]{bh1}.

Orbits of the projective elements are the most interesting from the
viewpoint of number theory in this case, and the transitivity result
stated in the theorem follows from Lemma~$\ref{red2}$ (it was also treated
in more detail in Theorem~$\ref{thm:main}$). The case $k=0$ corresponds to
the projective elements of rank~$3$ ($n=q'(x)=0$), and the case $k\ne 0$
corresponds to non-degenerate orbits ($q'(x)\ne 0$).

We notice that in the case $\jz=\h\binz$ the result of the theorem remains
valid for the subgroup $\grz\cc$ isomorphic to $\SL_6(\zet)$.
\ep

\section{Appendix: Bhargava's Cube Law and the Freudenthal construction}

In this Appendix we are going to present another link between
the higher composition laws and the Freudenthal construction.

M.~Bhargava discovered higher composition laws by providing a new perspective
on Gauss's composition, and then generalizing his results.
In order to define Gauss's composition, he considers three
quadratic forms associated to a cube with vertices labeled by integers.
He defines their addition (the Cube Law), and proves that this operation
is equivalent to Gauss's composition \cite{bh-phd, bh1}.

We also consider cubes of integers. Such cubes may be thought of as
elements of $\zet^2\tensor \zet^2 \tensor \zet^2$, and we use the
isomorphism of Example~$\ref{ex:FFF}$ to write such cubes as elements of
the module $\m(\jz)$ for $\jz=\zet\oplus\zet\oplus\zet$. Then it turns out
that the three quadratic forms of M.~Bhargava may be written using the
expressions that appeared in our definition of quadratic rank polynomials
(see Remarks~$\ref{rem:rp}$ and~$\ref{rem:rpz}$).

In a certain sense, Bhargava's construction is based on the three symmetries of
the cube with respect to the planes, parallel to the faces of the cube.
And our construction is based on the $3$-fold rotational
symmetry of the cube relative to one of its main diagonals.

\subsection{The original Cube Law}

The material from this section appeared in \cite[Section~2]{bh1}, see also
\cite[Section~2.1]{bh-phd}.

We consider the free $\zet$-module $\zet^2 \tensor \zet^2 \tensor \zet^2$
of rank $8$. This space has a natural basis, and an arbitrary element may
be written as an integral linear combination of the eight basis vectors.

It is convenient to represent elements of
$\zet^2 \tensor \zet^2 \tensor \zet^2$ by cubes, which have integer labels
attached to their vertices in the following way:
$$
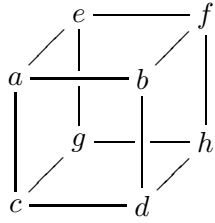
\xymatrix@!0{
  & e \ar@{-}[rr] \ar@{-}'[d][dd]
      &  & f \ar@{-}[dd]        \\
  a \ar@{-}[ur]\ar@{-}[rr]\ar@{-}[dd]
      &  & b \ar@{-}[ur]\ar@{-}[dd] \\
  & g \ar@{-}'[r][rr]
      &  & h                \\
  c \ar@{-}[rr]\ar@{-}[ur]
      &  & d \ar@{-}[ur]        }
$$

We consider three pairs $(M_i, N_i)$ of $2\times 2$ integer matrices
corresponding to the three possible slicings of the cube:
$$
M_1=\mat{a}bcd,\quad
N_1=\mat{e}fgh;
$$
$$
M_2=\mat{a}ceg,\quad
N_2=\mat{b}dfh;
$$
$$
M_3=\mat{a}ebf,\quad
N_3=\mat{c}gdh.
$$

We then construct three binary quadratic forms in the following way:
$$
Q_i(x,y)=-\det(M_i x -N_i y), \qquad
1\le i \le 3.
$$

M.~Bhargava introduced an operation ``+" on the set of (primitive) binary
quadratic forms by the relation (the Cube Law):
\begin{equation} \label{cubelaw}
Q_1 + Q_2 +Q_3 = 0.
\end{equation}

He proved that this operation is equivalent to Gauss's Law of composition
of quadratic forms, and it turns the set of $SL_2(\zet)$-equivalent
classes of primitive quadratic forms of discriminant $D$ into an abelian group.
This group is isomorphic to the (narrow) class group of a quadratic order
of discriminant $D$.

\subsection{A new realization of the Cube Law}

We again start with cubes with vertices labeled by integers in the
following way:
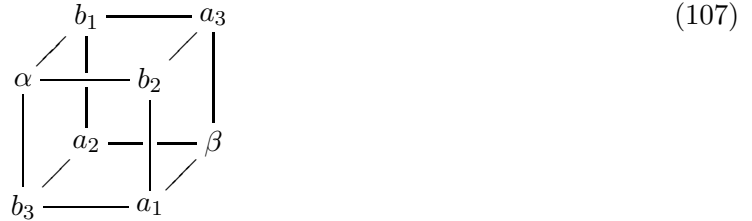
\begin{equation}\label{cube}
\xymatrix@!0{
  & b_1 \ar@{-}[rr] \ar@{-}'[d][dd]
      &  & a_3 \ar@{-}[dd]        \\
  \alpha \ar@{-}[ur]\ar@{-}[rr]\ar@{-}[dd]
      &  & b_2 \ar@{-}[ur]\ar@{-}[dd] \\
  & a_2 \ar@{-}'[r][rr]
      &  & \beta                \\
  b_3 \ar@{-}[rr]\ar@{-}[ur]
      &  & a_1 \ar@{-}[ur]        }
\end{equation}

Integers $\alpha$ and $\beta$ are located in the opposite vertices (they
correspond to the axis of rotation). We assign integers $a_i$ to the
vertices adjacent to $\beta$, and $b_i$'s
to the vertices adjacent to $\alpha$. Finally, we
choose labeling so that $a_i$ and $b_i$ are located at the opposite
vertices for $i=1,2,3$.

We will assign the following element $x\in \m(\jz),
\jz=\zet\oplus\zet\oplus\zet$
to the cube~$(\ref{cube})$:
$$
x=\mel\alpha\beta A B,
$$
where it is convenient to think of
$A$ and $B$ as being embedded diagonally
into the space of $3\times 3$ matrices
$$
A=
\matthree{a_1}\cdot  \cdot
         \cdot {a_2} \cdot
         \cdot \cdot {a_3}
\qquad
B=
\matthree{b_1}\cdot  \cdot
         \cdot {b_2} \cdot
         \cdot \cdot {b_3}.
$$

We consider the expressions (cf.~$(\ref{p2})$)
$$
\alpha A - B^\#, \quad
\beta  B - A^\#, \quad
Q(x),
$$
and note that for the diagonal matrices as above, $Q(x)$ becomes
the operator of multiplication by the diagonal matrix:
$$
(\alpha\beta-(A,B))I_3+ 2AB.
$$
Here $(A,B)=a_1b_1+a_2b_2+a_3b_3$ is the trace bilinear form, $AB$
represents the usual associative product of matrices,
and $A^\#$ is the usual adjoint: 
$A^\#=\diag{a_2a_3} {a_3a_1} {a_1a_2}$. 

\medskip
Consider the following expression:
\begin{equation} \label{main}
\biggl(\alpha A - B^\#\biggr)\,x^2 -
\biggl((\alpha\beta-(A,B))I_3+ 2AB\biggr) \,xy +
\biggl(\beta B - A^\#\biggr) y^2.
\end{equation}

The result is a $3\times 3$ diagonal matrix with entries:
$$
\matthree
{-R_1}  \cdot  \cdot
\cdot  {-R_2}  \cdot
\cdot  \cdot  {-R_3}
,
$$
where
\begin{eqnarray*}
-R_1 &=&
\biggl(\alpha a_1-b_2 b_3\biggr)x^2 +
\biggl(-a_1b_1+a_2 b_2+ a_3 b_3 -\alpha\beta\biggr)xy +
\biggl(\beta b_1-a_2a_3\biggr) y^2,\\[2mm]
-R_2 &=&
\biggl(\alpha a_2-b_3 b_1\biggr)x^2 +
\biggl(\phantom{-}a_1b_1-a_2 b_2+ a_3 b_3 -\alpha\beta\biggr)xy +
\biggl(\beta b_2-a_3a_1\biggr) y^2,\\[2mm]
-R_3 &=&
\biggl(\alpha a_3-b_1 b_2\biggr)x^2 +
\biggl(\phantom{-}a_1b_1+a_2 b_2- a_3 b_3 -\alpha\beta\biggr)xy +
\biggl(\beta b_3-a_1a_2\biggr) y^2.
\end{eqnarray*}

The quadratic forms $R_1, R_2, R_3$ are exactly the three forms $Q_1, Q_2,
Q_3$ appearing in the Cube Law.

{\sc
Department of Mathematics and Statistics,
University of Ottawa

585 King Edward Ave.,
Ottawa, ON, K1N 6N5, Canada
}

\medskip
E-mail: {\tt sergei.krutelevich@science.uottawa.ca}
\end{document}